\documentclass[10pt]{article}
\usepackage[utf8]{inputenc}
\usepackage{amssymb}
\usepackage{amsmath}
\usepackage{amsthm}
\usepackage{tikz}
\usepackage{hyperref}
\usetikzlibrary{arrows,calc,decorations.pathreplacing,fadings,intersections}
\usetikzlibrary{external}

\newtheorem{theorem}{Theorem}
\newtheorem{lemma}[theorem]{Lemma}
\newtheorem{corollary}[theorem]{Corollary}
\newtheorem{proposition}[theorem]{Proposition}
\newtheorem{definition}{Definition}

\newcommand{\vol}{\mathrm{vol}}

\title{ 
Ball-convex bodies and $L_p$ relative surface areas}
\author{
Elisabeth M. Werner \thanks{Supported by NSF grant DMS-2506790} \footnote{Department of Mathematics, Case Western Reserve University, 2145 Adalbert Road, Cleveland, OH 44106, USA,
{\tt elisabeth.werner@case.edu}}  \hskip 7mm
\hskip 7mm Diliya Yalikun \footnote{Department of Mathematics, Case Western Reserve University, 2145 Adalbert Road, Cleveland, OH 44106, USA, {\tt dxy259@case.edu@case.edu}} 
}
\date{}

\begin{document}

\maketitle

\begin{abstract}
We define new  surface area measures for ball-convex bodies which
we call  $L_p$ relative  surface  areas.
We show that those are  rigid motion invariant  valuations.
We establish inequalities  for these  quantities and prove a monotonicity behavior
which leads to a new notion of entropy for ball-convex bodies.
We introduce a weighted ball floating body. A derivative of volume of a ball-convex body with a  weighted ball floating body 
provides a geometric interpretation of the $L_p$ relative  surface  areas.

\end{abstract}

    

\section{Introduction} 

A convex body  in $\mathbb R^{n}$ is a compact, convex subset of $\mathbb R^{n}$ with nonempty
interior. Convex bodies $K$ in   $\mathbb R^{n}$  of the form
\begin{equation}\label{Rconvex}
K= 
\bigcap_{K \subset R \, B^n_2 +x} R \, B^n_2 + x, 
\end{equation}
where $R$ is a positive real number and $B^n_2$ is the Euclidean unit ball, are called $R$-ball convex. 
We denote the class of such bodies by $\mathcal{K}_R$.
 This  class has been introduced and investigated in \cite{BlaschkeKK, LangiNazodiTalata} and intensively studied, e.g., \cite{ArtsteinFlorentin, BezdekConnellyCsikos, BezdekNazodi, BezdekLangiNazodi, KabluchkoMarynychMolchanov, Pach}.
There are connections of this class to optimal transport, see \cite{ArtsteinSadovskyWyczensany}, to the Kneser-Poulsen conjecture (see \cite{Bezdek}), to Blaschke's rolling theorem  \cite{Drach}, to Meissner polyhedra \cite{Hynd1} and bodies of constant width  \cite{Hynd2, MartiniMontejanoOliveros}, and isoperimetric problems,   e.g., \cite{DrachTatarko}.
A classification of the isometries on ball-convex bodies  was given in  \cite{ArtsteinChorFlorentin1, ArtsteinChorFlorentin2}. 
In particular,  in \cite{SWY} a notion of a floating body has been introduced for this class, namely, the {\em $R$-ball floating body} of $K$, 
\begin{equation}\label{R-float}
K_\delta^R = \bigcap_{\vol_n(K\setminus R\, B^n_2 + x) \leq \delta} R\, B^n_2 +x , 
\end{equation}
where the Euclidean balls $R\, B^n_2+x$ replace the hyperplanes $H$ of the definition of the (classical) floating body, 
$K_\delta= \bigcap _{\vol_n(H^- \cap K) \leq \delta} H^+$, \cite{BaranyLarman1988, SW:1990}. $H^+$ and $H^-$ are  the two closed halfspaces
determined by the hyperplane $H$.
Note that when $R\to \infty$, we recover the classical floating body.
It was also shown in \cite{SWY} that a right derivative of volume of an $R$-ball convex body with its $R$-ball floating body exists,  
\begin{equation}\label{R-asa}
\lim_{\delta \to 0} \frac{\vol_n\left(K\right)- \vol_n\left( K_\delta^R \right)} {\delta^\frac{2}{n+1}} =  \alpha_n \int_{\partial K}  \prod _{i=1}^{n-1} \left( \kappa_i (K, x) -\frac{1}{R} \right)^\frac{1}{n+1} d \mu_K(x), 
\end{equation}
where $\alpha_n$ is a constant depending on $n$ only, $\mu_K$ is the usual surface area measure on the boundary $\partial K$ of $K$ and $\kappa_i(K, x)$, $1 \leq i \leq n-1$,  are the principal curvatures  of $K$ at $x$.
Their product is the Gauss curvature $\kappa (K,x)$.
The corresponding  result with   the classical floating body was proved in \cite{SW:1990} and gives as the limit the classical affine surface area which goes back to Blaschke \cite{Blaschke:1923}, 
\begin{equation}\label{asa}
as(K) = \int_{\partial K} \kappa(K, x)^\frac{1}{n+1} d\mu_K(x).
\end{equation}
As $R \to \infty$, the right hand side of (\ref{R-asa}) goes to (\ref{asa}) which motivated us to call the right hand side of (\ref{R-asa}) {\em relative affine surface area}.  Note though that it is not an affine invariant.
\par
\noindent
Due to its important properties, which make it an effective and powerful tool, the classical affine surface area 
is omnipresent in geometry, e.g., 
\cite{GHSW:2020, Hoehner:2022, LR:2010, Lutwak:1996, LYZ:2000, SW:2018}
and there are numerous  applications for it, 
such as, 
approximation  of convex bodies by polytopes, e.g., 
\cite{BLW:2018, BGT, Boe1, 
Reitz1, Reitzner, SW4}, Minkowski problems \cite{GuoXiZhao, HuangLutwakYangZhang}, affine curvature flows, information theory, and differential equations \cite{Andrews:1999, BesauWerner, TW1, TW2, TW4}.
The new relative surface area measure (\ref{R-asa}) has  also already appeared in approximation questions \cite{FodorKeveiVigh, FodorPapvari, FodorGrunfelder}.
\par
\noindent
An extension of the classical Brunn Minkowski theory,  the   {\em $L_p$ Brunn Minkowski theory}, was initiated by Lutwak in the groundbreaking paper \cite{Lutwak:1996}. At the core  of the $L_p$ Brunn Minkowski theory are the $L_p$ affine surface areas,
see below for the definition. They were introduced by Lutwak for $p >1$  in \cite{Lutwak:1996}, for $p>0$ in \cite{Hug}  and 
 in \cite{SW:2004} for all $p$. The case $p=1$ is the (classical) affine surface area (\ref{asa}).
The $L_p$ Brunn Minkowski theory  is now  a central part of modern convex geometry \cite{HaddadLangharstPuttermanRoysdonYe, LangharstXi,  SchuettWerner2023, TW:2019, TW:2023, WY:2008, Ye:2015, Ye:2016}.
Extensions of both, floating body and $L_p$ affine surface area, to spherical and hyperbolic space were developed in \cite{BW:2016, BW:2018, 
BesauWerner2024/2}.
\vskip 2mm
\noindent
The $L_p$ analogs for the relative surface  measures have been missing so far. Therefore it is natural to ask for $L_p$ versions in the context of $R$-ball convex bodies.
We introduce them in this paper.
\vskip 3mm
\noindent
The paper is organized as follows. In  section \ref{results} we present (some of) our results. In section \ref{Lp-def} we give the definition of  the $L_p$-relative  surface areas, 
show some of their properties and prove inequalities for these quantities. In section \ref{entropy} we introduce a new notion of entropy for ball-convex bodies. In section \ref{geometric}
we define weighted ball floating bodies and describe their relevance for  the $L_p$-relative  surface areas. In section \ref{proofs} we present the remaining proofs.
\vskip 3mm
\noindent
{\bf Further notation.} 
The  closed Euclidean ball centered at $a$ with radius $r$ is $B^n_2(a,r)$. We write in short $B^n_2=B^n_2(0,1)$ and $S^{n-1}= \partial B^n_2$. The Euclidean norm on $\mathbb{R}^n$ is $\| \cdot\|$.
We denote by  $H\left(x, \xi \right)$ the hyperplane through $x$ orthogonal to the vector $\xi$. 
The line segment joining $x$ and $y$ is $[x,y]$.
For a convex body $K$ in $\mathbb{R}^n$, $\vol_n(K)$, or $|K|$ in short, denotes it's volume. 
By $N(x)$ or $N_K(x)$ we denote the unit outer normal to $K$ in the boundary point $x$.
We will  assume throughout the paper without loss of generality that $0 \in \text{int}(K)$, the interior of $K$. The polar or dual of $K$ is
$$
K^\circ=\{y \in \mathbb{R}^n: \langle x, y \rangle \leq 1\}.
$$
For more information and details on convex bodies we refer to e.g., the books  \cite{Gardner,  SchneiderBook}.
\par\noindent
Finally,  $c$, $d$, $d_1$, $d_2$, are absolute constants that may change from line to line.
\vskip 3mm
\noindent
\section{Results} \label{results}
We now mention some of the results of the paper.  
\par
\noindent 
We denote by 
$\mathcal{K}_R^+$ the set of all $R$-ball convex bodies in $\mathbb{R}^n$
such that for all $x \in \partial K$, $\kappa_i(K,x) >\frac{1}{R}$, for all $1 \leq i \leq n-1$. 
\vskip 2mm
\noindent
The first theorem leads to our definition of $L_p$ relative surface areas.
We apply a method that was successful not only in the classical case \cite{MW:2000}, but also in spherical and hyperbolic space \cite{BesauWerner1, BesauWerner},  
namely we compose the floating  body operation with polarity. 
\vskip 3mm
\begin{theorem}\label{theorem:dual limit} 
Let $K \in \mathcal{K}_R^+$ and let $K_\delta^R$ be its $R$-ball floating body. Then
$$
\lim_{\delta \to 0} \frac{\vol_n\left((K_\delta^R)^\circ\right)- vol_n\left((K)^\circ\right)} {\delta^\frac{2}{n+1}} =  c_n \int_{\partial K} \frac{\kappa(K, x)}{\langle x, N(x)\rangle^{n+1}} \prod _{i=1}^{n-1} \left( \kappa_i (K, x) -\frac{1}{R} \right)^\frac{1}{n+1} d \mu_K(x), 
$$
where $c_n=\frac{n^2(n^2-1)^\frac{2}{n+1}}{2\, (\sigma (S^{n-2}))^\frac{2}{n+1}}$.
\end{theorem}
\noindent
Note that if $R \to \infty$, we obtain on the right hand side the $-n/(n+2)$-affine surface area \cite{Lutwak:1996, MW:2000, SW:2004}.
This leads naturally to calling  the integral expression the  {\em $-n/(n+2)$-relative surface area} and thus to define for all $p\neq -n$ 
the {\em $L_p$-relative  surface areas} of $R$-ball convex bodies $K$ as  
\begin{equation*}
\Omega_p^R(K) = \int_{\partial K}\Big (\frac{\kappa(K, x)^\frac{1}{n+1}}{\langle x, N(x)\rangle}\Big )^\frac{n(p-1)}{n+p}  \prod _{i=1}^{n-1} \left( \kappa_i (K, x) -\frac{1}{R} \right)^\frac{1}{n+1} d \mu_K(x).
\end{equation*}
For  $R \to \infty$, we recover the $L_p$-affine surface areas, defined by Lutwak \cite{Lutwak:1996} (see also \cite{Hug, SW:2004}), 
\begin{equation*}
as_p(K) = \int_{\partial K}\frac{\kappa(K, x)^\frac{p}{n+p}}{\langle x, N(x)\rangle ^\frac{n(p-1)}{n+p} } d \mu_K(x).
\end{equation*}
\vskip 3mm
\noindent
We show inequalities  for the $L_p$ relative surface areas, 
among them are  those in Theorem \ref{Ungleichungen}.
\vskip 2mm
\noindent
{\bf Theorem 4.} 
{\em Let $s> -n$, $r> -n$, $t> -n$, be real numbers. Let $K$ be an $R$-ball convex body in $\mathbb R^n$ such that  $0 \in \text{int} (K)$ 
and $\mu_K(\{x \in \partial K: \kappa_i(K,x)>\frac{1}{R}, \,  \text{for  all }  i,  1 \leq i \leq n-1 \})>0$.
\vskip 1mm
\noindent
(i) If $1 < \frac{(n+r)(t-s)}{(n+t)(r-s)} < \infty$, then
\begin{equation*}
    \Omega_r^R(K)\leq (\Omega_t^R(K))^\frac{(r-s)(n+t)}{(t-s)(n+r)} (\Omega_s^R(K))^\frac{(t-r)(n+s)}{(t-s)(n+r)}.
\end{equation*}
\vskip 1mm
\noindent
(ii) If $\frac{(n+r)t}{(n+t)r}>1$, then 
\begin{equation*}
    \left(\frac{\Omega^R_r(K)}{\Omega_0^R(K)}\right)\leq\left(\frac{\Omega^R_t(K)}{\Omega_0^R(K)}\right)^\frac{r(n+t)}{t(n+r)}.
\end{equation*}
\vskip 2mm
\noindent
Equality holds in both inequalities if and only if $K$ is an ellipsoid.
}
\vskip 3mm
\noindent
From the inequalities of Theorem \ref{Ungleichungen} we deduce a monotonicity  behavior  in $p$ of  $\left(\frac{\Omega_p^R(K)}{\Omega^R_\infty(K)}\right)^{n+p}$
which leads to the existence of the limit
$$
{E}^R(K)=\lim_{p\rightarrow\infty}\left(\frac{\Omega^R_p(K)}{\Omega^R_\infty(K)}\right)^{n+p}.
$$
${E}^R(K)$ is an entropy power and related to the Kullback Leibler divergence $D_{KL}(P||Q)$ of certain probability measures $P$ and $Q$,  see Section \ref{entropy}. This is shown in
Proposition \ref{KL}.
\vskip 2mm
\noindent
{\bf Proposition 8.}
{\em Let $K$ be a  $R$-ball convex body such that $0 \in \text{int} (K)$ 
    and  $\mu_K(\{x \in \partial K: \kappa_i(K,x)>\frac{1}{R}, \,  \text{for  all }  i,  1 \leq i \leq n-1 \})>0$. Then
\begin{equation*}\label{prop3:eq1} 
D_{KL}(P\|Q) 
 = \log\left( \frac{\Omega^R_0(K)} {\Omega^R_\infty(K)} \, {E}^R(K)^{-\frac{1}{n}}\right).
\end{equation*}
}
\vskip 2mm
\noindent
Finally, in Definition \ref{WFB}, we introduce {\em weighted $R$-ball  floating bodies} $F_R(K,f,\delta)$ where a set of volume $\delta$ weighted with a function $f$ is cut from an $R$-ball convex body.
Volume derivatives then lead to a weighted relative surface area of Theorem \ref{limit2}.  As a corollary (Corollary \ref{cor}) we obtain a geometric characterization of the $L_p$ relative surface areas.
\vskip 2mm
\noindent
{\bf Theorem 10.}
{\em Let $K$ be  a $R$-ball convex body in $\mathbb R^n$. Let $f:K\rightarrow \mathbb R$ be an continuous  function such that $f\geq c$ on $K$ where $c>0$ is a constant. Then
    $$\lim_{\delta\rightarrow 0}\frac{\vol_n(K)-\vol_n(F_R(K,f,\delta))}{\delta^{\frac{2}{n+1}}}= c_n\int_{\partial K}(f(x))^{-\frac{2}{n+1}}\prod_{i=1}^{n-1}\left(\kappa_i(K,x)-\frac{1}{R}\right)^\frac{1}{n+1}d\mu_{K}(x),$$
where $c_n=\frac{1}{2} \left( \frac{n+1}{\vol_{n-1}(B^{n-1}_{2})}\right)^\frac{2}{n+1}$.}

\vskip 5mm

\section{$L_p$-relative  surface areas} \label{Lp-def}

\begin{definition}\label{def:p-relative asa}
Let $K \in \mathcal{K}_R$.  Let $ - \infty \leq p \leq \infty$, $p\neq -n$. Then the  $L_p$-relative  surface area of $K$ is 
\begin{equation*}
\Omega_p^R(K) = \int_{\partial K}\Big (\frac{\kappa(K, x)^\frac{1}{n+1}}{\langle x, N(x)\rangle}\Big )^\frac{n(p-1)}{n+p}  \prod _{i=1}^{n-1} \left( \kappa_i (K, x) -\frac{1}{R} \right)^\frac{1}{n+1} d \mu_K(x).
\end{equation*}
\end{definition}
\vskip 2mm
\noindent
When $K$ is sufficiently smooth, $\Omega_p^R$ can be written as
\begin{equation} \label{Sn}
\Omega_p^R(K) = \int_{S^{n-1}} \prod _{i=1}^{n-1} \left( \frac{1}{r_i (K, u)} -\frac{1}{R} \right)^\frac{1}{n+1} \frac{f_K(u)^\frac{n^2+2n+p}{(n+1)(n+p)}}{h_K(u)^\frac{n(p-1)}{n+p}} d \sigma(u),
\end{equation}
where, for $u \in S^{n-1}$ such that $N_K(x)=u$,  $r_i (K, u)=(\kappa_i(K, N_K^{-1}(u))^{-1}$,  $f_K= \kappa^{-1}$ is the curvature function and $\sigma$ is the surface
area measure of the sphere.
\vskip 3mm
\noindent
{\bf Remarks and Examples.}
\vskip 3mm
\noindent
1. For $p=1$, we obtain the {\em relative ``affine" surface area} $as^R(K)$ introduced in \cite{SWY},
\begin{eqnarray*}\label{p=1}
as^R(K)=\Omega_1^R(K) &=& \int_{\partial K}  \prod _{i=1}^{n-1} \left( \kappa_i (K, x) -\frac{1}{R} \right)^\frac{1}{n+1} d \mu_K(x).
\end{eqnarray*}
\vskip 2mm
\noindent
2. If  $R >1$, $B^n_2$ is $R$-ball convex and then we have for all $p \neq -n$
\begin{equation}\label{ball}
\Omega_p^R(B^n_2)=\Omega_1^R(B^n_2).
\end{equation}
Next, we look at the case of an ellipsoid $\mathcal{E}$  in $\mathcal{K}_R$. We first recall a theorem by Petty \cite{Petty}.

\begin{theorem} \cite{Petty} \label{Petty} 
Let $K$ be a convex body in $C^2_+$. $K$ is an ellipsoid if and only if  for all $x$  in $\partial K$
$$
\frac{\kappa(K, x))^\frac{1}{n+1}} {\langle x, N(x) \rangle  } = k,$$
where $k>0$ is a constant.
\end{theorem}
\noindent
Thus we get for all $p \neq -n$
\begin{equation}\label{ellipse1}
\Omega_p^R(\mathcal{E})= k^\frac{n(p-1)}{n+p} \, \Omega_1^R(\mathcal{E}), 
\end{equation}
and in particular
\begin{equation}\label{ellipse2}
\Omega_{\pm \infty}^R(\mathcal{E})= k^{n} \, \Omega_1^R(\mathcal{E}), \hskip 12mm \Omega_{0 }^R(\mathcal{E})= k^{-1} \, \Omega_1^R(\mathcal{E}).
\end{equation}
\vskip 2mm
\noindent
3. Note that 
\begin{equation}\label{weightedp}
\Omega_p^R(K) = \int_{\partial K}\frac{\kappa(K, x)^\frac{p}{n+p}}{\langle x, N(x)\rangle^\frac{n(p-1)}{n+p} }\prod _{i=1}^{n-1} \left( 1 -\frac{1}{R\, \kappa_i (K, x)} \right)^\frac{1}{n+1} d \mu_K(x).
\end{equation}
For $R\to \infty$ we recover the usual $L_p$-affine surface area \cite{Lutwak:1996, Hug, SW:2004}.
\par
\noindent
In particular,  for $p=1$, $p=0$ and $p= \pm \infty$, we obtain 
\begin{eqnarray}\label{p=1}
as^R(K)=\Omega_1^R(K) 
&=&\int_{\partial K} \kappa (K, x)^\frac{1}{n+1} \prod _{i=1}^{n-1} \left( 1 -\frac{1}{R\, \kappa_i (K, x)} \right)^\frac{1}{n+1} d \mu_K(x),
\end{eqnarray}
\begin{eqnarray}\label{p=0}
\Omega_0^R(K) &=& 
\int_{\partial K} \langle x, N(x)\, \rangle \prod _{i=1}^{n-1} \left( 1 -\frac{1}{R\, \kappa_i (K, x)} \right)^\frac{1}{n+1} d \mu_K(x)
\end{eqnarray}
and  
\begin{eqnarray}\label{p=infty}
\Omega_{\pm \infty}^R(K)  
&=&\int_{\partial K}\frac{\kappa(K,x)}{\langle x, N(x)\, \rangle^n}\,  \prod _{i=1}^{n-1} \left( 1 -\frac{1}{R\, \kappa_i (K, x)} \right)^\frac{1}{n+1}  d \mu_K(x).
\end{eqnarray}
Therefore $\Omega_p^R(K)$ can be viewed as a weighted $L_p$-``affine" surface 
area, with the weight
\begin{equation}\label{weight}
w_K(x)=\prod _{i=1}^{n-1} \left( 1 -\frac{1}{R\, \kappa_i (K, x)} \right)^\frac{1}{n+1}.
\end{equation}
In particular,  as $n\, \vol_n(K) = \int_{\partial K} \langle x, N(x)\, \rangle d \mu_K(x)$, (\ref{p=0}) can be viewed as a weighted volume of $K$ and as 
 $n\, \vol_n(K^\circ) = \int_{\partial K} \frac{\kappa(K,x)}{\langle x, N(x)\, \rangle^n}\,  d \mu_K(x)$, (\ref{p=infty}) can be viewed as a weighted volume of $K^\circ$.
\vskip 2mm
\noindent
4. Let $K$ be a $R$-ball convex body. Then $\kappa_i(K,x) \geq \frac{1}{R}$ for all $i$, for all $x$.
It can happen for an $R$-ball convex body that one or all of the $\kappa_i(x) =\frac{1}{R}$.
If one or all of the $\kappa_i(x)=\frac{1}{R}$ almost everywhere on $\partial K$, then $\Omega_p^R(K)=0$, for all $p$. 
For instance this is the case for {\em $R$-ball polyhedra},  i.e. the intersection of
finitely many $R$-balls. In particular $\Omega_p^R(R B^n_2)=0$ for all $p$.
\vskip 2mm
\noindent

\vskip3mm

\subsection{Boundedness Property of the $L_p$-relative surface areas}

It was shown in \cite{Lutwak:1996} ($p > 1$; $p=1$ is the classical case going back to Blaschke \cite{Blaschke:1923}), in \cite{Hug} ($0 <p <1$) and for all  $p$
in \cite{WY:2008} that for a convex body $K$ with $0 \in \text{int}(K)$, 
\begin{eqnarray}\label{p-inequalities1}
\text{If}\, \,  p\geq 0: \, \frac{as_p(K)}{as_p(B^n_2)}\leq \left(\frac{|K|}{|B^n_2
|}\right)^{\frac{n-p}{n+p}}; \hskip 5mm \text{If}-n<p<0: \, \frac{as_p(K)}{as_p(B^n_2)}\geq
\left(\frac{|K|}{|B^n_2|}\right)^{\frac{n-p}{n+p}}.
\end{eqnarray}
Equality holds in both cases if and only if $K$ is an ellipsoid. 
For $p=0$, equality holds trivially for all $K$.
If  $K$ is in addition  in $C^2_+$ and if $p < -n$, then for an absolute constant $\gamma>0$, 
\begin{equation*}\label{p-inequalities2}
\gamma^{\frac{np}{n+p}}\left(\frac{|K|}{|B^n_2|
}\right)^{\frac{n-p}{n+p}} \leq \frac{as_p(K)}{as_p(B^n_2 )}. 
\end{equation*}
\vskip 2mm
\noindent
For an $R$-ball convex body we have that $\kappa_i (K,x) \geq \frac{1}{R}$ for all $i$, for all $x$.
Therefore, $\kappa(K,x) \geq \frac{1}{R^{n-1}}$ for all $x \in \partial K$. Moreover 
$\langle x, N(x)\rangle \leq R$ and thus, together with (\ref{p-inequalities1}),  we get that  
\begin{eqnarray*}
0 \leq \left(\frac{1}{R}\right)^{\frac{n^2(p-1)}{(n+1)(n+p)}} as_1^R(K) \leq \Omega_p^R(K) \leq & n |B^n_2|^\frac{2p}{n+p} |K|^ \frac{n-p}{n+p}, \hskip 5mm 1\leq p \leq \infty  
\end{eqnarray*}
\begin{eqnarray*}
0 \leq  \Omega_p^R(K) \leq & n |B^n_2|^\frac{2p}{n+p} |K|^ \frac{n-p}{n+p}, \hskip 5mm 0\leq p \leq 1
 \end{eqnarray*}
\begin{equation*}
0 \leq \Omega_p^R(K) \leq \left(\frac{1}{R}\right)^{n\frac{(n-1)(p-1)}{n+1)(n+p)}} 
\Omega_1^R(K), \hskip 5mm  -n <p \leq 1.
\end{equation*}
Thus the $\Omega_p^R(K) < \infty$ for $-n < p \leq \infty$.  In the range $ -\infty \leq p < -n$ we get
\begin{eqnarray*}
0 \leq \left(\frac{1}{R}\right)^{\frac{n^2(p-1)}{n+1)(n+p)}} \Omega_1^R(K) \leq \Omega_p^R(K) \leq as_p(K), \hskip 5mm -\infty \leq p < -n.
\end{eqnarray*}
It then can happen that $\Omega_p^R(K)$ is unbounded.  An example is $B^2_r =\{x \in \mathbb{R}^2: x_1^r +x_2^r \leq 1\}$, for $1 <r <2$. Then $B^2_r$ is $R$-ball   convex for $R\geq \frac{2^{-\frac{2-r}{2r}}}{r-1}$ but $\Omega_p^R(B^2_r)$ is unbounded when  $p(r-1)\geq -2$, and similarly for higher dimensions.
 \vskip 3mm
\noindent

\subsection{Homogeneity and Valuation Property
}

It is obvious that $L_p$-relative affine surface area is invariant under rigid motions. The next proposition states some more properties. 
\par
\begin{proposition}\label{rel-asa-properties} 
Let $K$ be an $R$-ball convex body and let $-\infty \leq p \leq \infty$, $p \neq -n$.
\vskip 1mm
\noindent
(i) Let $a \in \mathbb{R}$, $a >0$.  Then $\Omega_p^R$ is homogeneous of degree $n \frac{n-p}{n+p}$, 
\begin{equation*}
\Omega_p^{a R} (a K) = a^{n \frac{n-p}{n+p}} \, \Omega_p^R(K).
\end{equation*}
\vskip 1mm
\noindent
(ii) 
For $0 \leq p \leq \infty$, we have the following isoperimetric inequalities 
\begin{eqnarray*}
\Omega_p^R(K) \leq as_p(K) \leq n \, \vol_n(B^n_2)^\frac{2p}{n+p} \vol_n(K)^\frac{n-p}{n+p}.
\end{eqnarray*}
Equality holds in the first inequality iff $R=\infty$ and in the second inequality  iff $K$ is an ellipsoid. 
\vskip 1mm
\noindent
(iii) The $L_p$-relative affine surface area $\Omega_p^R$ is a valuation, i.e. for $R$-ball convex bodies $K$ and $L$ such that $K \cap L$ and $K\cup L$ are  again $R$-ball convex, 
$$
\Omega_p^R(K\cup L) + \Omega_p^R(K\cap L) = \Omega_p^R(K) +\Omega_p^R(L).
$$
\end{proposition}
\vskip 4mm
\subsection {Inequalities for $L_p$ -relative surface areas}
\vskip 2mm

\begin{theorem} \label{Ungleichungen}
    Let $s> -n$, $r>-n$, $t> -n$, be real numbers. Let $K$ be an $R$-ball convex body in $\mathbb R^n$ such that  $0 \in \text{int} (K)$ 
    and  $\mu_K(\{x \in \partial K: \kappa_i(K,x)>\frac{1}{R}, \,  \text{for  all }  i,  1 \leq i \leq n-1 \})>0$. 
\vskip 1mm
\noindent
(i) If $1 < \frac{(n+r)(t-s)}{(n+t)(r-s)} < \infty$, then
\begin{equation*}
    \Omega_r^R(K)\leq (\Omega_t^R(K))^\frac{(r-s)(n+t)}{(t-s)(n+r)} (\Omega_s^R(K))^\frac{(t-r)(n+s)}{(t-s)(n+r)}.
\end{equation*}
\vskip 1mm
\noindent
(ii) If $\frac{(n+r)t}{(n+t)r}>1$, then 
\begin{equation*}
    \left(\frac{\Omega^R_r(K)}{\Omega_0^R(K)}\right)\leq\left(\frac{\Omega^R_t(K)}{\Omega_0^R(K)}\right)^\frac{r(n+t)}{t(n+r)}.
\end{equation*}
\vskip 2mm
\noindent
Equality holds in both inequalities if and only if $K$ is an ellipsoid.
\end{theorem}
\vskip 2mm
\noindent
{\bf Remarks.}  
\par
\noindent
If $\mu_K(\{x \in \partial K: \kappa_i(K,x)>\frac{1}{R}, \,  \text{for \, all }  i,  1 \leq i \leq n-1 \})=0$, then $\Omega_p^R(K)=0$ for all $p$.
\par
\noindent
$\frac{(n+r)(t-s)}{(n+t)(r-s)}=1$ happens if either $s=-n$ which is excluded, or when $r=t$. In that case, equality holds trivially in (i).
\par
\noindent
Similarly, $\frac{(n+r)(t-s)}{(n+t)(r-s)}=\infty$ happens if either $t=-n$ which is excluded, or when $r=s$. In that case, again equality holds trivially in (i).
\par
\noindent
In the same way, $\frac{(n+r)t}{(n+t)r}=1$ means that $r=t$ and then equality holds trivially in (ii) and $\frac{(n+r)t}{(n+t)r}=\infty$, means that either $t=-n$ which is excluded,  or $r=0$ and then equality holds trivially in (ii).
\vskip 2mm
\begin{proof} We follow the proof given in  \cite{WY:2008} for the $L_p$ affine surface areas.
\vskip 2mm
\noindent
(i)  By
H\"older's inequality - which enforces the condition $\frac{(n+r)(s-t)}{(n+t)(s-r)} >1$
 \begin{eqnarray*} \label{firstcase}
 && \hskip -12mm \Omega_r^R(K)=\int _{\partial 
K} \Big (\frac{\kappa(K, x)^\frac{1}{n+1}}{\langle x, N(x)\rangle}\Big )^\frac{n(r-1)}{n+r}  \prod _{i=1}^{n-1} \left( \kappa_i (K, x) -\frac{1}{R} \right)^\frac{1}{n+1} d \mu_K(x) 
\\
&& \hskip -4mm = \int _{\partial K}
\left(\frac{\kappa(K, x)^{\frac{1}{n+1}}}{\langle x, N_K(x)\rangle}
\right)^{n\frac{(t-1)(r-s)(n+t)}{(n+t)(t-s)(n+r)}}
\left(\prod _{i=1}^{n-1} \left( \kappa_i (K, x) -\frac{1}{R} \right)^\frac{1}{n+1} \right)^\frac{(r-s)(n+t)}{(t-s)(n+r)} \\
&&\times \left(\frac{\kappa(K, x)^{\frac{1}{n+1}}}{\langle x, N_K(x)\rangle
}\right)^{n\frac{(s-1)(t-r)(n+s)}{(n+s)(t-s)(n+r)}}\,
\left(\prod _{i=1}^{n-1} \left( \kappa_i (K, x) -\frac{1}{R} \right)^\frac{1}{n+1}\right)^\frac{(t-r)(n+s)}{(t-s)(n+r)}
d\mu _K(x) \\
 && \hskip -4mm\leq  \big(\Omega^R_t(K)
\big)^{\frac{(r-s)(n+t)}{(t-s)(n+r)}} \big(\Omega^R_s
(K)\big)^{\frac{(t-r)(n+s)}{(t-s)(n+r)}}. 
 \end{eqnarray*}
 We have applied H\"older's inequality with $p=\frac{(t-s)(n+r)}{(r-s)(n+t)}$, $q=\frac{(t-s)(n+r)}{(t-r)(n+s)}$, 
 $$
f= \left(\frac{\kappa(K, x)^{\frac{1}{n+1}}}{\langle x, N_K(x)\rangle
} \right)^{n\frac{(t-1)(r-s)(n+t)}{(n+t)(t-s)(n+r)}}
 \left(\prod _{i=1}^{n-1} \left( \kappa_i (K, x) -\frac{1}{R} \right)^\frac{1}{n+1}
 \right)^\frac{(r-s)(n+t)}{(t-s)(n+r)}
 $$
$$
g=\left(\frac{\kappa(K, x)^{\frac{1}{n+1}}}{\langle x, N_K(x)\rangle
} \right)^{n\frac{(s-1)(t-r)(n+s)}{(n+s)(t-s)(n+r)}}
 \left(\prod _{i=1}^{n-1} \left( \kappa_i (K, x) -\frac{1}{R} \right)^\frac{1}{n+1}
 \right)^\frac{(t-r)(n+s)}{(t-s)(n+r)}.
$$
 Equality  holds in  H\"older's inequality if and only if  for some constant $c>0$, $|f|^p=c|g|^q$ a.e.
 i.e., if and only if 
 $$
 \frac{\kappa(K, x)^\frac{1}{n+1}}{\langle x, N_K(x)\rangle
}  = c^\frac{(n+s)(n+t)}{n(n+1)(t-s)}.
 $$   
By Petty's Theorem \ref{Petty} this happens if and only  if $K$ is an ellipsoid.
\vskip 2mm
\noindent
(ii) Similarly, again using H\"older's inequality - which now
enforces the condition $\frac{(n+r)t}{(n+t)r} >1,$
\begin{eqnarray*}
  \Omega_r^R(K)&=&\int _{\partial 
K} \Big (\frac{\kappa(K, x)^\frac{1}{n+1}}{\langle x, N(x)\rangle}\Big )^\frac{n(r-1)}{n+r}  \prod _{i=1}^{n-1} \left( \kappa_i (K, x) -\frac{1}{R} \right)^\frac{1}{n+1} d \mu_K(x) 
\\
&=& \int _{\partial K}
\left(\frac{\kappa(K, x)^{\frac{1}{n+1}}}{\langle x, N_K(x)\rangle
} \right)^{n\frac{(t-1)r(n+t)}{(n+t) t(n+r)}}
 \left(\prod _{i=1}^{n-1} \left( \kappa_i (K, x) -\frac{1}{R} \right)^\frac{1}{n+1}
 \right)^\frac{r(n+t)}{t(n+r)} \\
&&\hskip 5mm \times\left(\frac{\kappa(K, x)^{\frac{1}{n+1}}}{\langle x, N_K(x)\rangle
} \right)^{-n\frac{t-r}{t(n+r)}}
 \left(\prod _{i=1}^{n-1} \left( \kappa_i (K, x) -\frac{1}{R} \right)^\frac{1}{n+1}
 \right)^\frac{n(t-r)}{t(n+r)} 
 d\mu_K(x) \\
&\leq &\big(\Omega^R_t(K)\big)^{\frac{r(n+t)}{t(n+r)}} \big(\Omega_0^R(K)
\big)^{\frac{(t-r)n}{(n+r)t}}.
\end{eqnarray*}
 We have now applied H\"older's inequality with $p=\frac{t(n+r)}{r(n+t)}$, $q=\frac{t(n+r)}{n(t-r)}$, 
 $$
f= \left(\frac{\kappa(K, x)^{\frac{1}{n+1}}}{\langle x, N_K(x)\rangle
} \right)^{n\frac{(t-1)r(n+t)}{(n+t) t(n+r)}}
 \left(\prod _{i=1}^{n-1} \left( \kappa_i (K, x) -\frac{1}{R} \right)^\frac{1}{n+1}
 \right)^\frac{r(n+t)}{t(n+r)}
 $$
$$
g=\left(\frac{\kappa(K, x)^{\frac{1}{n+1}}}{\langle x, N_K(x)\rangle
} \right)^{-n\frac{t-r}{t(n+r)}}
 \left(\prod _{i=1}^{n-1} \left( \kappa_i (K, x) -\frac{1}{R} \right)^\frac{1}{n+1}
 \right)^\frac{n(t-r)}{t(n+r)}.
$$
Equality  holds in  H\"older's inequality if and only if 
for some constant $c>0$,
 $$
 \frac{\kappa(K, x)^\frac{1}{n+1}}{\langle x, N_K(x)\rangle
}  = c^\frac{n+t}{t(n+1)}.
 $$   
Again, by Petty's Theorem \ref{Petty},  this happens if and only  if $K$ is an ellipsoid.
\end{proof}
\vskip 1mm
\vskip 4mm
\section{Entropies and their inequalities} \label{entropy}
\vskip 2mm
It follows from Theorem \ref{Ungleichungen} (ii) that if
 $0<r<t$, or $-n < r <t < 0 $,  then
\begin{equation*}
    \left(\frac{\Omega^R_r(K)}{\Omega_0^R(K)}\right)^\frac{n+r}{r}\leq\left(\frac{\Omega^R_t(K)}{\Omega_0^R(K)}\right)^\frac{n+t}{t}.
\end{equation*}
It follows from Theorem \ref{Ungleichungen} (i) that if $s \to \infty$ and $\frac{n+r}{n+t}>1$, then 
$$\frac{\Omega_r^R(K)}{\Omega_\infty^R(K)}\leq\left(\frac{\Omega^R_t(K)}{\Omega_\infty^R(K)}\right)^\frac{n+t}{n+r}.$$
Thus the next proposition is a consequence of Theorem \ref{Ungleichungen}.
\vskip 3mm
\begin{proposition} \label{monoton}
    Let $K$ be an $R$-ball convex body in $\mathbb R^n$ with $0$ in it's interior
    and such that $\mu_K(\{x \in \partial K: \kappa_i(K,x)>\frac{1}{R}, \,  \text{for  all }  i,  1 \leq i \leq n-1 \})>0$. Then
\vskip 1mm
\noindent
(i) The function $p\rightarrow\left(\frac{\Omega_p^R(K)}{\Omega^R_\infty(K)}\right)^{n+p}$ is decreasing in $p\in(-n,\infty)$.
\vskip 1mm
\noindent
(ii) The function $p\rightarrow\left(\frac{\Omega_p^R(K)}{\Omega_0^R(K)}\right)^\frac{n+p}{p}$ is increasing in $p\in(-n,\infty)$.
\end{proposition}

\vskip 2mm
\noindent
  Petty's Theorem \ref{Petty}  implies that we have strict monotonicity in Proposition \ref{monoton} (i), (ii), unless $K$ is an ellipsoid, in which case the quantities (i) and (ii) are constant. By (\ref{ellipse1}) and (\ref{ellipse2})  the quantity  in  (i) equals $k^{-n(n+1)}$ and the quantity (ii) equals to $k^{p\frac{n+1}{n+p}}$.
\vskip 4mm
\noindent
Proposition \ref{monoton} leads to the next definition.

\vskip 2mm

\begin{definition} \label{def2}
    Let $K$ be an $R$-ball convex body in $\mathbb R^n$ with $0$ in it's interior
    and such that $\mu_K(\{x \in \partial K: \kappa_i(K,x)>\frac{1}{R}, \,  \text{for  all }  i,  1 \leq i \leq n-1 \})>0$. We define
    $$E^R(K)=\lim_{p\rightarrow\infty}\left(\frac{\Omega^R_p(K)}{\Omega^R_\infty(K)}\right)^{n+p}$$.
\end{definition}
\vskip 3mm
\noindent
 Note that by (\ref{ball}), (\ref{ellipse1}) and (\ref{ellipse2}) we have
that $E^R(B^n_2) =1$, when $R>1$ and that for an $R$-ball ellipsoid $\mathcal E$ that $E^R(\mathcal E) =\frac{1}{k^{n(n+1)}}$, where $k$ is the constant from 
Theorem \ref{Petty}.
\vskip 3mm
\begin{proposition} \label{decrease}
Let $K$ be an $R$-ball convex body in $\mathbb R^n$ with $0$ in it's interior and such that $\mu_K(\{x \in \partial K: \kappa_i(K,x)>\frac{1}{R}, \,  \text{for  all }  i,  1 \leq i \leq n-1 \})>0$. Then
\begin{eqnarray*}
&&\log E ^R(K)\\
&&=
 - \frac{n}{\Omega^R_{\infty}(K)} \int_{\partial K} 
\frac{\kappa(K, x)}{\langle x, N_{K}(x) \rangle^{n}} \log\left( \frac{\kappa(K, x)}{\langle x, N_{K}(x) \rangle^{n+1}}\right)
 \prod _{i=1}^{n-1} \left(1 -\frac{1}{R \, \kappa_i (K, x)} \right)^\frac{1}{n+1} d\mu_{K}(x).
\end{eqnarray*}
\end{proposition}
\vskip 3mm
\begin{proof}
We have by Definition \ref{def2} and de l'Hopital's rule that 
\begin{eqnarray*}
\hskip -5mm\log{E^R(K)} 
 &=& \log{\left( \lim_{p\rightarrow \infty}  \left(\frac{\Omega^R_{p}(K)}{\Omega^R_{\infty}(K)}\right)^{n+p} \right) }= \lim_{p \rightarrow \infty} \frac{ \log{ \left( \frac{\Omega^R_{p}(K)}{\Omega^R_{\infty}(K)}\right) } }{ (n+p)^{-1}} 
 \\
 &=&- \lim_{p \rightarrow \infty} \frac{ (n+p)^{2} \frac{d}{dp}\left(\Omega^R_{p}(K)\right) }{\Omega^R_{p}(K)} \\
&=& 
 - \lim_{p \rightarrow \infty} \frac{(n+p)^{2}}{\Omega^R_{p}(K)} \int_{\partial K} \frac{d}{dp} \left[\exp  \left(\frac{n(p-1)}{n+p} \log \left(\frac{\kappa(K,x)^\frac{1}{n+1}}{\langle x, N_K(x)\rangle}\right)\right) \right] \\
 &&\hskip 20mm \times\prod _{i=1}^{n-1} \left( \kappa_i (K, x) -\frac{1}{R} \right)^\frac{1}{n+1} 
d\mu_{K}(x)  \\
 &=&
 -\lim_{p \rightarrow \infty} \frac{ (n+p)^{2} n(n+1)}{(n+p)^{2}\Omega^R_{p}(K)} \int_{\partial K} \left(\frac{\kappa(K, x)^{\frac{1}{n+1}}}{\langle x, N_{K}(x) \rangle}\right)^{\frac{n(1-\frac{1}{p})}{1+\frac{n}{p}}} \log\left( \frac{\kappa(K, x)^\frac{1}{n+1}}{\langle x, N_{K}(x) \rangle}\right) \\
 & &\hskip 20mm \times \prod _{i=1}^{n-1} \left( \kappa_i (K, x) -\frac{1}{R} \right)^\frac{1}{n+1} d\mu_{K}(x)  \\
 &= &
 - \frac{n(n+1)}{\Omega^R_{\infty}(K)} \int_{\partial K} 
 \left(\frac{\kappa(K, x)^{\frac{1}{n+1}}}{\langle x, N_{K}(x) \rangle}\right)^{n} \log\left( \frac{\kappa(K, x)^\frac{1}{n+1}}{\langle x, N_{K}(x) \rangle}\right) \\
 & &\hskip 20mm \times \prod _{i=1}^{n-1} \left( \kappa_i (K, x) -\frac{1}{R} \right)^\frac{1}{n+1} d\mu_{K}(x)\\
&=& - \frac{n}{\Omega^R_{\infty}(K)} \int_{\partial K} 
\frac{\kappa(K, x)}{\langle x, N_{K}(x) \rangle^{n}} \log\left( \frac{\kappa(K, x)}{\langle x, N_{K}(x) \rangle^{n+1}}\right) \\
 & &\hskip 20mm \times \prod _{i=1}^{n-1} \left(1 -\frac{1}{R \, \kappa_i (K, x)} \right)^\frac{1}{n+1} d\mu_{K}(x).
 \end{eqnarray*}
\end{proof}
\noindent
For $R \to \infty$, we recover the entropy expression of \cite{PW}.
\vskip 3mm
\noindent
The next corollary follows immediately from Proposition \ref{monoton} and the definition of $E^R(K)$.
\begin{corollary} \label{cor-info} 
Under the same assumptions as in Proposition \ref{monoton}
we have for all $p \geq 0$
$$
E^R(K) \leq \left(\frac{\Omega_p^R(K)}{\Omega^R_\infty(K)}\right)^{n+p} \leq \left(\frac{\Omega_0^R(K)}{\Omega^R_\infty(K)}\right)^{n}, 
$$
with equality in the inequalities if and only if $K$ is an ellipsoid.
\end{corollary}
\vskip 4mm
\noindent 
Let $(X, \mu)$ be a measure space  and let  $dP=pd\mu$ and  $dQ=qd\mu$ be probability measures on $X$ that are  absolutely continuous with respect to the measure $\mu$. 
The {\em Kullback-Leibler divergence} or {\em relative entropy} from $P$ to $Q$ is defined as \cite{CT}
\begin{equation}\label{relent}
 D_{KL}(P\|Q)= \int_{X} p\log{\frac{p}{q}} d\mu.
\end{equation}
\vskip 3mm
\noindent 
{\em The information inequality} (also  called {\em  Gibb's inequality}) \cite{CT} holds for the Kullback-Leibler divergence:
Let $P$ and $Q$ be  as above. Then
\begin{equation}\label{Gibbs}
  D_{KL}(P\|Q) \ge 0,
\end{equation}
with equality if and only if $P=Q$.
\vskip 2mm
\noindent
Let $K$ be an $R$-ball convex  body in $\mathbb R^n$ such that $0 \in \text{int} (K)$ 
    and  $\mu_K(\{x \in \partial K: \kappa_i(K,x)>\frac{1}{R}, \,  \text{for  all }  i,  1 \leq i \leq n-1 \})>0$.
Let 
\begin{equation}\label{PQ}
p(x)= \frac{ \kappa (K, x)}{ \Omega_\infty^R(K)\,\langle x, N_{K}(x) \rangle^{n} \  } \, , \   \ q(x)= \frac{\langle x, N_{K}(x) \rangle}{ \Omega_0^R(K)}.
\end{equation}
Let 
$$
w_K(x)=\prod _{i=1}^{n-1} \left( 1 -\frac{1}{R\, \kappa_i (K, x)} \right)^\frac{1}{n+1}
$$
be the weight introduced in (\ref{weight}) and put $d\,  w_K= w_K d \mu_K$. 
Then $dP=p\  d \, w_K$ and $dQ=q \ d \, w_K$ are probability measures on $\partial K$ that are absolutely continuous with respect  to $\mu_{K}$.
\vskip 3mm
\noindent
The next proposition shows that $E^R$ is an entropy power.

\begin{proposition}\label{KL}
Let $K$ be an  $R$-ball convex body such that $0 \in \text{int} (K)$ 
    and  $\mu_K(\{x \in \partial K: \kappa_i(K,x)>\frac{1}{R}, \,  \text{for  all }  i,  1 \leq i \leq n-1 \})>0$. Then
\begin{equation*}\label{prop3:eq1} 
D_{KL}(P\|Q) 
 = \log\left( \frac{\Omega_0(K)} {\Omega_\infty(K)} \, E^R(K)^{-\frac{1}{n}} \right).
\end{equation*}
\end{proposition} 
\vskip 2mm
\noindent
By the information inequality we see, as in Corollary \ref{cor-info}, that 
$E^R(K) \leq \left(\frac{\Omega_0^R(K)}{\Omega^R_\infty(K)}\right)^{n}$, with equality if and only if $K$ is an ellipsoid.

\vskip 4mm
\section{Geometric Interpretations} \label{geometric}
\subsection{Weighted $C$-ball floating bodies}
\vskip 3mm
\begin{definition} \label{WFB} Let $C$ be a fixed convex body in $\mathbb R^n$ and
    let $K$ be a $C$-ball convex body. Let $0\leq \delta$ and let $f:K\rightarrow \mathbb R$ be an integrable strictly positive function. 
    We define the weighted $C$-ball floating body $F_C(K,f,\delta)$ of $K$ by
    $$F_C(K,f, \delta)=\bigcap_{\int_{K\setminus (C+x)} f\,dm\leq \delta}C+x.
    $$
    When $C=RB^n_2$, we write in short $F_R(K,f,\delta)=F_{RB^n_2}(K,f,\delta)$ and call it the weighted $R$-ball convex floating body.
\end{definition}
\vskip 3mm
\noindent
It is obvious that $F_C(K,f,0)=K$. Note also that for $f=1$ we recover the $R$-ball floating body of \cite{SWY}. It follows immediately from the transformation formula that for every  linear transformation $T$ with $\det T\neq 0$,
$$
F_{TC}(TK,f\circ T^{-1},\delta)=T\left(F_C\left(K,f,\frac{\delta}{|\det T|}\right)\right).
$$
\par
\noindent
The next proposition describes another property of the weighted $R$-ball floating bodies. The proof is the same as for the non-weighted version in \cite{SWY} but we include it for completeness. 
\par
\noindent
\begin{proposition}\label{properties}
Let $K$ be an $R$-ball convex body in $\mathbb R^n$   and let $f:K\rightarrow \mathbb R$ be a continuous  strictly positive function. 
     For all $\delta$ such that $F_R(K,f,\delta)\neq\emptyset$ and all $x_\delta\in\partial F_R(K,f,\delta) \cap\text{int}(K)$ there exists at least one $R$-ball $R\, B^n_2$ such that $\int_{K\setminus R\, B^n_2}f dm=\delta$.
\end{proposition}     
\vskip 2mm
\begin{proof}
    \vskip 1mm
    \noindent
    Without loss of generality we choose $\delta>0$ is so small that $F_R(K,f,\delta)\neq \emptyset$. In particular, we choose $\delta < \frac{1}{2}\int _K f dm$. 
Let $x\in \partial F_R(K,f,\delta)$. We choose a sequence $x_{k}$, $k\in \mathbb N$,
such that $x_{k}\notin F_R(K,f,\delta)$ and $\lim_{k\to\infty}x_{k}=x$. By the
definition of $F_R(K,f,\delta)$ we find for every $x_{k}$ a ball $B^n_2(y_k, R)$  such
that $x_{k}\notin B^n_2(y_k, R)$ and $\int_{K\setminus B^n_2(y_k, R))} f dm =\delta$. By
compactness there is a subsequence $B^n_2(y_{k_{j}}, R)$, $j \in \mathbb N$, that converges (in the Hausdorff metric) to
an $R$-ball $B^n_2(y,R)$, i.e. the sequence  of the centers converges to the center and there is a sequence of elements on $\partial B^n_2(y_{k_{j}}, R)$  that converges to
an element on  $\partial B^n_2(y, R)$.
Therefore $\int_{K \setminus B^n_2(y, R))}f dm = \delta$ and $x \in  \partial B^n_2(y, R)$.
\end{proof}
\vskip 3mm
\noindent

\begin{theorem} \label{limit2}
   Let $K$ be  a $R$-ball convex body in $\mathbb R^n$. Let $f:K\rightarrow \mathbb R$ be an continuous  function such that $f\geq c$ on $K$ where $c>0$ is a constant. Then
    $$\lim_{\delta\rightarrow 0}\frac{\vol_n(K)-\vol_n(F_R(K,f,\delta))}{\delta^{\frac{2}{n+1}}}= c_n\int_{\partial K}(f(x))^{-\frac{2}{n+1}}\prod_{i=1}^{n-1}\left(\kappa_i(K,x)-\frac{1}{R}\right)^\frac{1}{n+1}d\mu_{K}(x),$$
where $c_n=\frac{1}{2} \left( \frac{n+1}{\vol_{n-1}(B^{n-1}_{2})}\right)^\frac{2}{n+1}$.
\end{theorem}
\vskip 3mm
\noindent
\begin{corollary} \label{cor}
Let $K$ be  a $R$-ball convex body in $\mathbb R^n$. Let $p \in \mathbb R$, $p \neq -n$ and let 
 $f_p:K\rightarrow \mathbb R$ be such that for all $x \in \partial K$
$$
f(x) = \left(\frac{\langle x, N_K(x)\rangle ^{n(n+1)(p-1)}}{\kappa(K,x)^{n(p-1)}}\right)^\frac{1}{2(n+p)}.
$$
 Then, with $c_n$ as above, 
    $$\lim_{\delta\rightarrow 0}\frac{\vol_n(K)-\vol_n(F_R(K,f_p,\delta))}{\delta^{\frac{2}{n+1}}}= c_n\,  \Omega_p^R(K).
    $$
\end{corollary}

\vskip 4mm
\section{The Remaining Proofs} \label{proofs}
We will need more notations and lemmas and the following Proposition \ref{integral}. For a proof of Proposition \ref{integral} see e.g., \cite{SWY}.
\vskip 3mm
\begin {proposition} \label{integral}
For $1 \leq i \leq n$, let $c_i >0$. Then
$$
\int_{S^{n-1}} \frac{d\sigma(\xi)}
{\left(\sum_{i=1}^{n} c_i \xi_i^2\right)^\frac{n}{2}} = \frac{2 \pi^\frac{n}{2}}{\Gamma(\frac{n}{2})} \left[\prod_{i=1}^{n} c_i ^\frac{1}{2}\right]^{-1} = \sigma (S^{n-1}) \left[\prod_{i=1}^{n} c_i ^\frac{1}{2}\right]^{-1}.
$$
\end{proposition}
\vskip 3mm
\noindent
For  the next lemma, we recall the \emph{rolling function}, also called {\em interior reach} $r_K:\partial K \to [0,\infty)$ of a convex body $K$, which  was introduced by McMullen in \cite{MM:1974}, see also \cite{SW:1990}.
For $x \in \partial K$ with unique outer normal $N_K(x)$ it  is defined by 
\begin{equation*}
	r_K(x) = \max\{\rho: B^{n}_2 (x - \rho N_{K} (x), \rho ) \subset K\}, 
\end{equation*}
i.e., $r_K(x)$ is the maximal radius of a Euclidean ball inside $K$ that contains $x$.
If $N_K(x)$ is not unique, $r_K(x) =0$. By McMullen \cite{MM:1974} (also \cite{SW:1990}) $r_K(x)>0$  almost everywhere on $\partial K$.
It was shown in \cite{SW:1990} that for all $0 \leq \alpha <1$, 
\begin{equation}\label{r}
\int_{\partial K} \frac{1}{r_K(x)^\alpha} d \mu_K(x) < \infty.
\end{equation}
\vskip 4mm
\noindent
Throughout this section, $f:K\rightarrow \mathbb R$ be an continuous  function such that 
\begin{equation}\label{constantc}
f\geq c 
\end{equation}
on $K$ where $c>0$ is a constant.
\vskip 3mm
\noindent
The next lemma and its proof is the analog of Lemma 6 of \cite{SW:1990}, see also Lemma 6 of \cite{SWY}. 
\vskip 3mm
\begin{lemma} \label{bounded}
Let $K$ be an $R$-ball convex body in $\mathbb R^{n}$  that contains $0$ as an interior point. 
Let $x \in \partial K$ such that $r_{K}(x)>0$ and let $x_{\delta}$ be the unique point in $[0,x]\cap\partial F_R(K,f,\delta)$. Then there is $\delta_0$ such that for all $\delta \leq \delta_0$,  
\begin{equation}\label{bounded-1}
\frac{\|x-x_\delta\|}{\delta^\frac{2}{n+1}} \leq \gamma_n \, r_{K}(x)^{-\frac{n-1}{n+1}},
\end{equation}
where $\gamma_n$ depends on $n$ and $K$ only.
\end{lemma}
\vskip 2mm
\begin{proof}
Since $0$ is an interior point of $K$ there is $s>0$ such that 
\begin{equation}\label{s-contain}
B_{2}^{n}(0,\tfrac{1}{s})\subseteq K\subseteq B_{2}^{n}(0,s).
\end{equation}
We choose $\delta_{0}=\frac{c}{2}\text{vol}_{n}(B_{2}^{n}(0,\tfrac{1}{4s}))$, where $c$ is the constant from (\ref{constantc}).
In order to prove \eqref{bounded-1} we consider two cases,
 $s^{-2}r_{K}(x)\leq \|x-x_{\delta}\|_{2}$ and $s^{-2}r_{K}(x)\geq \|x-x_{\delta}\|$.
\par
\noindent
We assume first that $s^{-2}r_{K}(x)\leq \|x-x_{\delta}\|$.
More specifically, we will show 
\begin{equation}\label{FloatTh3-2}
\frac{\|x-x_\delta\|}{\delta^\frac{2}{n+1}}
\leq \frac{2^{2\frac{n-1}{n+1}} \, (n/c)^{\frac{2}{n+1}}s^{6\frac{n-1}{n+1}}}{(\text{vol}_{n-1}
(B_{2}^{n-1}))^{\frac{2}{n+1}}} \, r_{K}(x)^{-\frac{n-1}{n+1}}.
\end{equation}
Let $\delta \leq \delta_0$ and let
$z +R B^n_2$ be an $R$-ball that touches $F_R(K,f,\delta)$ in $x_\delta$ and cuts off a set of weighted volume exactly 
equal to $\delta$ from $K$,
i.e. $\int_{K\setminus (z +R B^n_2)} f dm =\delta$.
Such an $R$-ball exists by Proposition \ref{properties} (ii). 
If $z +R B^n_2$ does not contain $B_{2}^{n}(0,\tfrac{1}{2s})$ then
 $K\setminus (z +R B^n_2)$ contains a Euclidean ball of radius
$\frac{1}{4s}$.  Then, as $f \geq c>0$ on $K$, 
$$
\delta = \int_{K\setminus (z +R B^n_2)} f dm \geq  c\, \text{vol}_{n}(B_{2}^{n}(0,\tfrac{1}{4s})).
$$
This cannot be since $\delta$ is smaller than
$\delta_{0}=\frac{c}{2}\text{vol}_{n}(B_{2}^{n}(0,\tfrac{1}{4s}))$. Therefore
\begin{equation}\label{bounded-2}
 B_{2}^{n}(0,\tfrac{1}{2s})\subseteq z +R B^n_2.
\end{equation}
Let $T(x_\delta)$ be  the tangent hyperplane to $z +R B^n_2$ in $x_\delta$ and $T(x_\delta)^+$ the half space containing
 $z +R B^n_2$.
Then
$$
C=[x,B_{2}^{n}(0,\tfrac{1}{2s})]\cap T(x_\delta)^{-}
$$
is a cone with a base that is an ellipsoid. Instead of the cone $C$
we consider the cone 
$$
\tilde C
=[x,B_{2}^{n}(0,\tfrac{1}{2s})]\cap H(x_{\delta})^-, 
$$
where $H(x_{\delta})^{}$ is the hyperplane through $x_\delta$ and orthogonal to $x$. By a simple geometric argument
$$
\operatorname{vol}_{n}(\tilde C)\leq\operatorname{vol}_{n}(C).
$$
 The height of the  cone $\tilde C$ is $\|x-x_{\delta}\|$ and the radius of the
base is $\frac{1}{2s}\frac{\|x-x_{\delta}\|}{\|x\|}$.
Therefore, using again that $f \geq c>0$ on $K$, 
\begin{equation*}\label{FloatTh3-1}
 \delta  = \int_{K\setminus (z +R B^n_2)} f dm  \geq c\, \operatorname{vol}_{n}(K\cap  H(x_{\delta})^-)
\geq \frac{c}{n}\frac{\|x-x_{\delta}\|^{n}}{\|x\|^{n-1}}\operatorname{vol}_{n-1}
(B_{2}^{n-1}(0,\tfrac{1}{2s})).
\end{equation*}
And thus
$$
\frac{\|x-x_\delta\|}{\delta^\frac{2}{n+1}} 
\leq \frac{(n/c)^{\frac{2}{n+1}} \|x\|^{2\frac{n-1}{n+1}}(2s)^{2\frac{n-1}{n+1}}}{(\text{vol}_{n-1}
(B_{2}^{n-1}))^{\frac{2}{n+1}}\|x-x_{\delta}\|^{\frac{n-1}{n+1}}}.
$$
Since $\|x\| \leq s$ by (\ref{s-contain}) and since
$s^{-2}r_{K}(x)\leq \|x-x_{\delta}\|$ by assumption, we get
\begin{eqnarray*}
\frac{\|x-x_\delta\|}{\delta^\frac{2}{n+1}}
\leq \frac{(n/c)^{\frac{2}{n+1}} 2^{2\frac{n-1}{n+1}}s^{4\frac{n-1}{n+1}}}{(\text{vol}_{n-1}
(B_{2}^{n-1}))^{\frac{2}{n+1}}\|x-x_{\delta}\|^{\frac{n-1}{n+1}}} 
\leq \frac{(n/c)^{\frac{2}{n+1}} 2^{2\frac{n-1}{n+1}}s^{6\frac{n-1}{n+1}}}{(\text{vol}_{n-1}
(B_{2}^{n-1}))^{\frac{2}{n+1}}} r_{K}(x)^{-\frac{n-1}{n+1}}.
\end{eqnarray*}
 This proves (\ref{FloatTh3-2}).
\vskip 1mm
\noindent
Now we 
consider the case $s^{-2}r_{K}(x)\geq \|x-x_{\delta}\|$.
Again, let $T(x_\delta)$ denote the tangent hyperplane to $z +R B^n_2$ in $x_\delta$. Then, since $B_{2}^{n}(x-r_{K}(x)N_K(x),r_{K}(x))\subseteq K$
\begin{eqnarray}\label{cap}
\delta&=&
 \int_{K\setminus (z +R B^n_2)} f dm \geq \int_{(K\cap  T(x_{\delta})^-)} f dm
\nonumber
 \geq \int_{(B_{2}^{n}(x-r_{K}(x)N_K(x),r_{K}(x))\cap T(x_{\delta})^-)} f dm \\
&\geq& c\, \vol_n\left((B_{2}^{n}(x-r_{K}(x)N_K(x),r_{K}(x))\cap T(x_{\delta})^-)\right).
\end{eqnarray}
Now we are in the situation of Lemma 6 of \cite{SW:1990} resp.  Lemma 6 of \cite{SWY} with $H=T(x_{\delta})$. For completeness, we give the arguments.
\par
\noindent
Let $\theta$ denote the angle between $x$ and $N_{K}(x)$. 
By \eqref{s-contain} and \eqref{bounded-2} we have 
\begin{equation}\label{bound-5}
\cos\theta\geq\frac{1}{s\|x\|}\geq\frac{1}{s^{2}}.
\end{equation}
We show that
\begin{equation}\label{bound-3}
x_{\delta}\in B_{2}^{n}(x-r_{K}(x)N_K(x),r_{K}(x)).
\end{equation}
We have (the inequality follows by \eqref{bound-5})
\begin{equation}\label{bound-8}
\operatorname{vol}_{1}\left([0,x]\cap B_{2}^{n}(x-r_{K}(x)N_K(x),r_{K}(x)-{K}(x))\right)
=2 r_{K}(x)\cos\theta\geq \frac{2 r_{K}(x)}{s^{2}}.
\end{equation}
 By the assumption
$s^{-2}r_{K}(x)\geq \|x-x_{\delta}\|$ we get \eqref{bound-3}.
\par
Let $\Delta$ be the distance from $x_{\delta}$ to the boundary of
$B_{2}^{n}(x-r_{K}(x)N_K(x),r_{K}(x))$
$$
\Delta =\min\{\|x_{\delta}-z\|\ : \ \|z-(x-r_{K}(x)N_K(x))\|=r_{K}(x)\}.
$$
$B_{2}^{n}(x-r_{K}(x)N_K(x),r_{K}(x))\cap T(x_{\delta})^-$ is a cap of height $\Delta$ of the ball $B_{2}^{n}(x-r_{K}(x)N_K(x),r_{K}(x))$.
We use  e.g., Lemma 8 of  \cite{SW:1990}  to estimate the volume of this cap  
\begin{eqnarray*}
\delta &\geq& c\, \vol_{n}(B_{2}^{n}(x-r_{K}(x)N_K(x),r_{K}(x))\cap T(x_{\delta})^-) \\
&\geq& 2c\, (2r_{K}(x))^{\frac{n-1}{2}}\frac{\operatorname{vol}_{n-1}(B_{2}^{n-1})}{n+1}
\left\{\Delta^{\frac{n+1}{2}}-\frac{(n+1)(n-1)}{4r_{K}(x) (n+3)}
\Delta^{\frac{n+3}{2}}\right\}
\\
&=&2c\, (2r_{K}(x))^{\frac{n-1}{2}}\frac{\operatorname{vol}_{n-1}(B_{2}^{n-1})}{n+1}
\Delta^{\frac{n+1}{2}}\left\{1-\frac{(n+1)(n-1)}{4r_{K}(x) (n+3)}
\Delta^{}\right\}.
\end{eqnarray*}
\begin{figure}[h]
			\centering
	\includegraphics[scale=0.5]{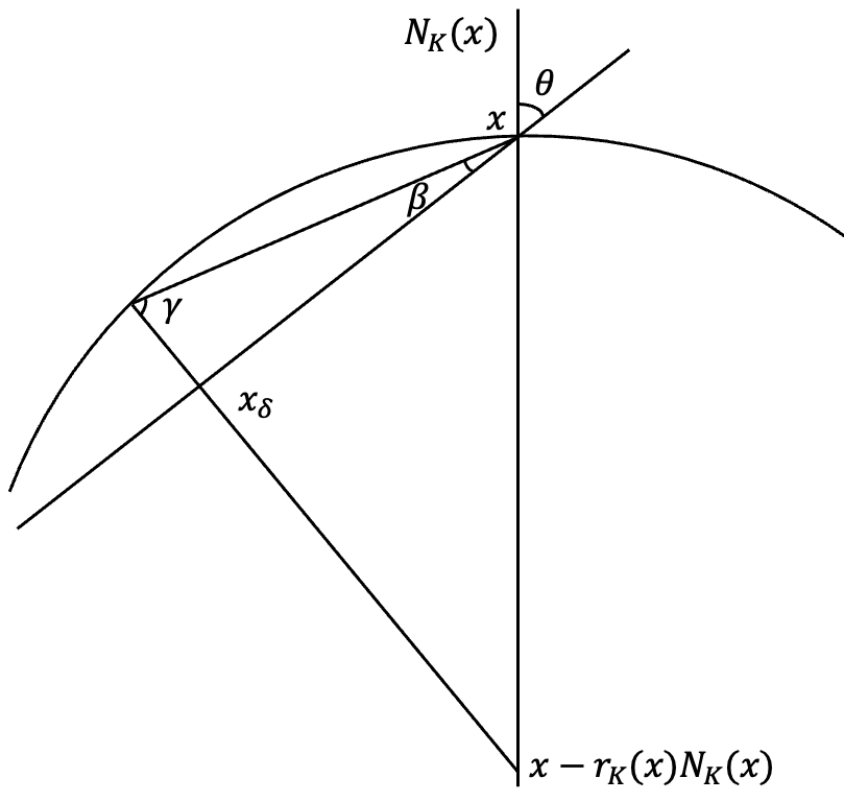}		
		\caption{The angles}
		\end{figure}	
Let $\theta$,  $\beta$  and $\gamma$ be as in Figure 1. 
If $\theta=0$ then $\Delta=\|x-x_{\delta}\|$. If $\Delta>0$, then $\gamma>0$.
By the sine law for the triangle
\begin{equation*}\label{FloatTh3-4}
\Delta=\frac{\sin \beta}{\sin\gamma}\|x-x_{\delta}\|.
\end{equation*}
Since $\beta\leq\gamma\leq\frac{\pi}{2}$ we have
$\sin\beta\leq\sin\gamma$. Therefore, $\Delta\leq\|x-x_{\delta}\|$. Moreover,
$\gamma\geq \theta$. Together with \eqref{s-contain} this implies $\sin \gamma\geq\sin\theta\geq\frac{1}{s^{2}}$. Thus
$$
\delta  \geq 
2c\,(2r_{K}(x))^{\frac{n-1}{2}}\frac{\operatorname{vol}_{n-1}(B_{2}^{n-1})}{n+1}
(\sin\beta \|x-x_{\delta}\|)^{\frac{n+1}{2}}\left\{1-\frac{(n+1)(n-1)}{4 \, r_{K}(x) (n+3)}
  \|x-x_{\delta}\|\right\}  .
$$
It follows that
$$
r_{K}(x)^{-\frac{n-1}{n+1}}  
\geq 2c\,  \left(\frac{\operatorname{vol}_{n-1}(B_{2}^{n-1})}{n+1}\right)^{\frac{2}{n+1}}
\frac{\sin\beta \|x-x_{\delta}\|_{2}}{ \delta^{\frac{2}{n+1}}}
\left\{1-\frac{(n+1)(n-1)}{4\, r_{K}(x) (n+3)}
  \|x-x_{\delta}\|\right\}^{\frac{2}{n+1}}.
$$
Since $s^{-2}r_{K}(x)\geq\|x-x_{\delta}\|$, 
\begin{eqnarray*}
r_{K}(x)^{-\frac{n-1}{n+1}}  
&\geq& 2c\,  \left(\frac{\vol_{n-1}(B_{2}^{n-1})}{n+1}\right)^{\frac{2}{n+1}}
\frac{\sin\beta \|x-x_{\delta}\|}{ \delta^{\frac{2}{n+1}}}
\left\{1-\frac{(n+1)(n-1)}{4\, s^{2} (n+3)}
\right\}^{\frac{2}{n+1}}, 
\end{eqnarray*}
which implies that, choosing $s \geq (n+1)^\frac{1}{2}$, 
\begin{eqnarray*}
\frac{ \|x-x_{\delta}\|}{ \delta^{\frac{2}{n+1}}} \leq \frac{1}{\sin\beta }\, \frac{1}{2c} \left(\frac{n+1}{\vol_{n-1}(B_{2}^{n-1})}\right)^{\frac{2}{n+1}} \left(\frac{4}{3}\right)^{-\frac{2}{n+1}}  r(x)^{-\frac{n-1}{n+1}}.
\end{eqnarray*}
It is left to verify that $\sin\beta$ is bounded from below. 
Let $\alpha$ be the angle between $x$ and $x_{\delta}$ at $x-r_{K}(x)N_K(x)$.
We have
\begin{equation}\label{bound-9}
\gamma=\beta+\theta
\hskip 20mm
2\gamma+\alpha=\pi
\hskip 20mm
\alpha+\theta+\frac{\pi}{2}\leq\pi
\end{equation}
Indeed, $\gamma=\beta+\theta$ follows immediately from Figure 1. For $2\gamma+\alpha=\pi$
we use that the angle sum in a triangle equals $\pi$. By \eqref{bound-8} and by the assumption
$s^{-2}r_{K}(x)\geq \|x-x_{\delta}\|$ we conclude that
the angle at $x_{\delta}$ between $x$ and $x-r_{K}(x)N_K(x)$ is greater than $\frac{\pi}{2}$.
This implies $\alpha+\theta+\frac{\pi}{2}\leq\pi$.
\par
By \eqref{bound-9}
$$
\tfrac{1}{2}(\tfrac{\pi}{2}-\theta)\leq\beta.
$$
On the interval $[0,\frac{\pi}{2}]$ the function $\sin$
is increasing. By \eqref{bound-5}
$$
\sin(2\beta)\geq\sin(\tfrac{\pi}{2}-\theta)
=\sin(\tfrac{\pi}{2})\cos(-\theta)
+\cos(\tfrac{\pi}{2})\sin(-\theta)
=\cos\theta\geq\tfrac{1}{s^{2}}.
$$
\end{proof}
\vskip 3mm
\noindent
\subsection{Proof of Theorem \ref{theorem:dual limit}}
\vskip 2mm
\noindent
We can assume that $0 \in \text{int}(K)$. Let $\delta >0$. 
We recall that $\vol_n(K^\circ) =\int _{S^{n-1}} \frac{1}{h_K^n(u)} d \sigma(u)$, where $h_K(u) = \sup_{ x \in K} \langle x, u \rangle$ is the support function of $K$. Therefore, 
\begin{equation}\label{voldifference}
\vol_n\left((K^R_\delta)^\circ\right)-\vol_n\left( K^\circ \right)=\int_{S^{n-1}}\left(\frac{1}{h_{K^R_\delta}^n(u)}-\frac{1}{h_{K}^n(u)}\right)d\sigma(u).
\end{equation}
For  $u \in S^{n-1}$, let 
$x \in \partial K$ be such that $N_K(x)=u$. Let $x_\delta \in \partial K_\delta^R$ be such that $x_\delta =[0,x] \cap \partial K_\delta^R$. 
Then, with an absolute constant $d$, 
\begin{equation}
n \, \|x-x_\delta\| \left\langle \frac{x}{\|x\|}, u \right\rangle (1 - d \|x-x_\delta\|)\leq h_K(u)-h_{K^R_\delta}(u) \leq n \, \|x-x_\delta\| \left\langle \frac{x}{\|x\|}, u \right\rangle.
\end{equation}
By Lemma 8 of \cite {SWY} and Proposition \ref{integral} we have that 
\begin{eqnarray*}
\lim_{\delta \to 0} \left\langle \frac{x}{\|x\|}, u \right\rangle  \frac{ \|x-x_\delta\|}{\delta^\frac{2}{n+1}} &=&  
\frac{\left(n^{2}-1\right)^\frac{2}{n+1} }{2}
\left[\int_{S^{n-2}} \frac{d\sigma(\xi)}
{\left(\sum_{i=1}^{n-1} \left(\kappa_i(x, K) -\frac{1}{R} \right) \xi_i^2\right)^\frac{n-1}{2}}\right]^{-\frac{2}{n+1}}\\
&=&\frac{\left(n^{2}-1\right)^\frac{2}{n+1} }{2\, (\sigma (S^{n-2}))^\frac{2}{n+1}} \prod_{i=1}^{n-1} \left(\kappa_i(x,K) - \frac{1}{R}\right)^\frac{1}{n+1}.
\end{eqnarray*}
Therefore, with a (new) constant $d$,
\begin{eqnarray*}
&&\hskip -15mm \frac{n(n^2-1)^\frac{2}{n+1}}{2\, (\sigma (S^{n-2}))^\frac{2}{n+1}} (1-d \, \delta^\frac{2}{n+1})  \delta^\frac{2}{n+1}   \prod_{i=1}^{n-1} \left(\kappa_i(x,K) - \frac{1}{R}\right)^\frac{1}{n+1}  
\leq h_K(u)-h_{K^R_\delta}(u)   \\
&& \hskip 35mm \leq \frac{n(n^2-1)^\frac{2}{n+1}}{2\, (\sigma (S^{n-2}))^\frac{2}{n+1}}   \delta^\frac{2}{n+1}   \prod_{i=1}^{n-1} \left(\kappa_i(x,K) - \frac{1}{R}\right)^\frac{1}{n+1}, 
\end{eqnarray*}
or
\begin{eqnarray*}
&&\frac{1}{h_{K}^n(u)}\left[ 1+\frac{n^2(n^2-1)^\frac{2}{n+1}}{2\, (\sigma (S^{n-2}))^\frac{2}{n+1}h_K(u)\,} (1-d_1 \, \delta^\frac{2}{n+1})  \delta^\frac{2}{n+1} \prod_{i=1}^{n-1} \left(\kappa_i(x,K) - \frac{1}{R}\right)^\frac{1}{n+1}\right]  \\
&&\leq \frac{1}{h_{K^R_\delta}^n(u)}\\
&&\leq \frac{1}{h_{K}^n(u)}\left[ 1+\frac{n^2(n^2-1)^\frac{2}{n+1}}{2\, \, (\sigma (S^{n-2}))^\frac{2}{n+1}h_K(u)} \left(1+d_2 \delta^\frac{2}{n+1}\right)  \delta^\frac{2}{n+1}   \prod_{i=1}^{n-1} \left(\kappa_i(x,K) - \frac{1}{R}\right)^\frac{1}{n+1} \right]. 
\end{eqnarray*}
By (\ref{voldifference}) and Lemma \ref{bounded} which we apply with $f=1$ on $K$, 
\begin{eqnarray*}
&&\lim_{\delta \to 0} \frac{\vol_n\left((K^R_\delta)^\circ\right)-\vol_n\left(K^\circ\right)}{\delta^\frac{2}{n+1}} = \int_{S^{n-1}}
\lim_{\delta \to 0} \frac{1}{\delta^\frac{2}{n+1}}\left(\frac{1}{h_{K^R_\delta}^n(u)}-\frac{1}{h_{K}^n(u)}\right)d\sigma(u)\\
&& = \lim_{\delta \to 0} \int_{S^{n-1}}
 \frac{1}{\delta^\frac{2}{n+1}}\left(\frac{1}{h_{K^R_\delta}^n(u)}-\frac{1}{h_{K}^n(u)}\right)d\sigma(u)\\
&&= \frac{n^2(n^2-1)^\frac{2}{n+1}}{2\, (\sigma (S^{n-2}))^\frac{2}{n+1}}   \int_{S^{n-1}} \frac{ \prod_{i=1}^{n-1} \left(\kappa_i(N_K^{-1}(u),K) - \frac{1}{R}\right)^\frac{1}{n+1}}{h_{K}^{n+1}(u)} d \sigma (u)\\
&&= \frac{n^2(n^2-1)^\frac{2}{n+1}}{2\, (\sigma (S^{n-2}))^\frac{2}{n+1}}   \int_{\partial K}\frac{\kappa(x, K)}{\langle x, N_K(x)\rangle^{n+1} }  \prod_{i=1}^{n-1} \left(\kappa_i(x,K) - \frac{1}{R}\right)^\frac{1}{n+1} d \mu_K (x).
\end{eqnarray*}

\vskip 4mm
\subsection{Proof of Proposition \ref{rel-asa-properties}}
\vskip 2mm
\noindent
{\bf Proof of (i)}
\par
\noindent
$K$ is $R$-ball convex is equivalent to $a K$ is $a R$-ball convex. 
We now use that $\mu_{a K}= a^{n-1} \mu_K$, that for  $y=a \, x$,  $x\in \partial K$, 
$\langle y, N_{aK}(y) \rangle = a\langle x, N_{K}(x) \rangle$, and that for $1 \leq i \leq n-1$
$$ 
 \kappa_i(aK,y) = \frac{1}{a} \kappa_i(K,x), \hskip 2mm \kappa(aK, y) = \frac{\kappa(K,x)}{a^{n-1}}.
$$
Then (i) follows immediately.
\vskip 3mm
\noindent
{\bf Proof of (ii)}
\par
\noindent
The first inequality is obvious. The second inequality is the  $L_p$- affine isoperimetric inequality 
and equality there holds iff $K$ is an ellipsoid \cite{Lutwak:1996, Hug, WY:2008}.
\vskip 3mm
\noindent
{\bf Proof of (iii)}
\par
\noindent
We follow the approach of Sch\"utt \cite{Sc} and decompose the involved boundaries into disjoint sets, namely 
\begin{eqnarray} \label{AffSurfVal11-15}
		\partial(C\cup K)
		&=&\{\partial C\cap\partial K\}\cup\{\partial C\cap K^{c}\}\cup\{C^{c}\cap \partial K\} \nonumber  \\
		\partial(C\cap K)
		&=&\{\partial C\cap\partial K\}
		\cup\{\partial C\cap\operatorname{int}(K)\}\cup\{\operatorname{int}(C)\cap \partial K\}   \nonumber \\
		\partial C&=&\{\partial C\cap\partial K\}\cup\{\partial C\cap K^{c}\}
		\cup\{\partial C\cap\operatorname{int}(K)\}  \nonumber \\
		\partial K&=&\{\partial C\cap\partial K\}\cup\{\partial K\cap C^{c}\}
		\cup\{\partial K\cap\operatorname{int}(C)\},
	\end{eqnarray}
see e.g.,  \cite{SWY} for a detailed proof.
Moreover, the sets $\partial K\cap\partial L$, $\partial K\cap L^{c}$ and $\partial  L\cap K^{c}$ are pairwise disjoint, and the sets
$\partial K\cap\partial L$, $\partial K\cap \operatorname{int}(L)$ and $\partial  L\cap \operatorname{int}(K)$ are pairwise disjoint.
We have for all  $i=1,\dots,n-1$,  for all $x\in\partial L\cap K^{c}$, respectively, for all $x\in \partial K\cap L^{c}$, 
\begin{equation} \label{hauptkruemmungen}
\kappa_i (K\cup L, x)=\kappa_i ( L, x), \hskip 2mm \text{respectively,} \hskip 2mm\kappa_i (K\cup L, x)=\kappa_i ( K, x).
\end{equation}
We have for all $i=1,\dots,n-1$, for almost all $x\in\partial K\cap\partial L$, 
\begin{equation}\label{AffSurfVal11-32}
\kappa_{i}(K\cup L,x)=\kappa_{i}(K,x)
\hskip10mm\mbox{and}\hskip10mm
\kappa_{i}(K\cap L,x)=\kappa_{i}(L,x)
\end{equation}
or  for  all $i=1,\dots,n-1$
\begin{equation}\label{AffSurfVal11-33}
\kappa_{i}(K\cup L,x)=\kappa_{i}(L,x)
\hskip10mm\mbox{and}\hskip10mm
\kappa_{i}(K\cap L,x)=\kappa_{i}(K,x).
\end{equation}
For a proof see e.g., \cite{SWY}
We also have, see e.g., \cite{Sc} for a proof,
for all $x\in\partial L\cap\partial K$ where all the curvatures $\kappa(L\cup K, x)$, $\kappa(L\cap K, x)$,
$\kappa(L,x)$ and $\kappa(K, x)$ exist that
\begin{equation}\label{kruemmung}
\kappa(L\cup K, x)=\min\{\kappa(L, x),\kappa(K, x)\}
\hskip 15mm
\kappa(L\cap K, x)=\max\{\kappa(L, x),\kappa(K, x)\}.
\end{equation}
We now write $\mu_{\partial K}$ instead of $\mu_K$ so that it is clear
which surface area measure is meant. 
By \eqref{AffSurfVal11-15}, \ref{hauptkruemmungen}, \ref{AffSurfVal11-32},  \ref{AffSurfVal11-33} and \ref{kruemmung}
\begin{eqnarray*}
&&\Omega_p^R(K\cup L)=
\int_{\partial (K\cup L)} \frac{\kappa(K\cup L, x)^\frac{n(p-1)}{(n+1)(n+p)}}{\langle x, N(x)\rangle^\frac{n(p-1)}{n+p} } \prod _{i=1}^{n-1} \left( \kappa_i (K\cup L, x) -\frac{1}{R} \right)^\frac{1}{n+1}
d\mu_{\partial (K\cup L)}\\
&&=
\int_{\partial K\cap\partial L}\frac{\min\left\{\kappa(L, x)^\frac{n(p-1)}{(n+1)(n+p)},\kappa(K, x)^\frac{n(p-1)}{(n+1)(n+p)}\right\}}{\langle x, N(x)\rangle^\frac{n(p-1)}{n+p} }\prod _{i=1}^{n-1} \left( \kappa_i (K\cup L, x) -\frac{1}{R} \right)^\frac{1}{n+1}
d\mu_{\partial (K\cup L)}
\\
&&
+\int_{\partial L\cap K^{c}}\frac{\kappa(L, x)^\frac{n(p-1)}{(n+1)(n+p)}}{\langle x, N(x)\rangle^\frac{n(p-1)}{n+p} }\prod _{i=1}^{n-1} \left( \kappa_i ( L, x) -\frac{1}{R} \right)^\frac{1}{n+1}
d\mu_{\partial L}\\
&&+\int_{\partial  K\cap L^{c}}\frac{\kappa(K, x)^\frac{n(p-1)}{(n+1)(n+p)}}{\langle x, N(x)\rangle^\frac{n(p-1)}{n+p} }\prod _{i=1}^{n-1} \left( \kappa_i (K, x) -\frac{1}{R} \right)^\frac{1}{n+1}
d\mu_{\partial K}.
\end{eqnarray*}
and 
\begin{eqnarray*}
&&\Omega_p^R(K\cap L) = \int_{\partial (K\cap L)}\frac{\kappa(K\cap L, x)^\frac{n(p-1)}{(n+1)(n+p)}}{\langle x, N(x)\rangle^\frac{n(p-1)}{n+p} }\prod _{i=1}^{n-1} \left( \kappa_i (K\cap L, x) -\frac{1}{R} \right)^\frac{1}{n+1}
d\mu_{\partial (K\cap L)}\\
&&=\int_{\partial K\cap\partial L}\frac{\max\left\{\kappa(L, x)^\frac{n(p-1)}{(n+1)(n+p)},\kappa(K, x)^\frac{n(p-1)}{(n+1)(n+p)}\right\}}{\langle x, N(x)\rangle^\frac{n(p-1)}{n+p} }\prod _{i=1}^{n-1} \left( \kappa_i (K\cap L, x) -\frac{1}{R} \right)^\frac{1}{n+1}
d\mu_{\partial (K\cap L)}
\\
&&
+\int_{\partial L\cap \operatorname{int}(K)}\frac{\kappa(L, x)^\frac{n(p-1)}{(n+1)(n+p)}}{\langle x, N(x)\rangle^\frac{n(p-1)}{n+p} }\prod _{i=1}^{n-1} \left( \kappa_i ( L, x) -\frac{1}{R} \right)^\frac{1}{n+1}d\mu_{\partial L}\\
&&
+\int_{\partial  K\cap \operatorname{int}(L)}\frac{\kappa(K, x)^\frac{n(p-1)}{(n+1)(n+p)}}{\langle x, N(x)\rangle^\frac{n(p-1)}{n+p} }\prod _{i=1}^{n-1} \left( \kappa_i ( K, x) -\frac{1}{R} \right)^\frac{1}{n+1}d\mu_{\partial K}.
\nonumber
\end{eqnarray*}
Note that on the set $\partial K\cap\partial L$ the measures $\mu_{\partial (K\cup L)}$, $\mu_{\partial (K\cap L)}$,
$\mu_{\partial K}$ and $\mu_{\partial L}$
all coincide.
Therefore, 
\begin{eqnarray*}
\Omega_p^R(K\cup L) + \Omega_p^R(K\cap L) = \Omega_p^R(K) +\Omega_p^R(L).
\end{eqnarray*}
We have also used that for any real numbers $a$ and $b$, 
	$$
	a+b=\min\{a,b\}+\max\{a,b\}.
	$$
\vskip 4mm
\subsection{Proof of Theorem \ref{limit2}}
\vskip 3mm
\noindent
The next lemma is standard, see e.g., \cite{SW:1990}.   

\begin {lemma} \label{vol-diff}
Let $K$ and $L$  be a convex bodies in $\mathbb{R}^n$ such that $0 \in \text{int} (L) \subset K$. Then
$$
\vol_n(K)-\vol_n(L) = \frac{1}{n} \int_{\partial K} \langle x, N_K(x) \rangle \left[1- \left(\frac{\|x_L\|}{\|x\|}\right)^n \right] d \mu_K(x), 
$$
where $x_L=[0,x] \cap \partial L$.
\end{lemma} 
\vskip 2mm
\noindent
We can assume without loss of generality that $0 \in \text{int} (F_R(K,f,\delta))$ and when $L=F_R(K,f,\delta)$, we write $x_\delta= x_{F_R(K,f,\delta)}$.
\vskip 3mm
\noindent
We say that $R\, B^n_2$ is a supporting $R$-ball to $K$ if $R\, B^n_2 \cap \partial K \neq \emptyset$ and $K \subset R\, B^n_2$.
\vskip 3mm
\noindent
In the next lemma we estimate the volume of a cap $C(h, \mathcal{E})$ of height $h$ cut off  from an ellipsoid by an $R$-ball.
More precisely, let an ellipsoid $\mathcal{E}$  be given by 
\begin{equation}\label{ellipse}
\sum_{i=1}^{n-1} \frac{x_i^2}{a_i^2} + \frac{(x_n-a_n)^2}{a_n^2}=1
\end{equation}
and let the $R$-ball have equation
\begin{equation}\label{cut-ball}
\sum_{i=1}^{n-1}x_i^2 + (x_n-a)^2=R^2,
\end{equation}
i.e., the $R$-ball is centered at $z=(0 \cdots, 0,a)$.
Without loss of generality, we can assume that $a > R \geq a_n \geq a_{n-1} \dots \geq a_1$.
Then  
$$
C(h, \mathcal{E}) = \mathcal{E} \setminus B^n_2(z,R)
$$
and the height $h$ of the cap is $h=a-R$. 
\vskip 2mm
\noindent
The next lemma and its proof is in \cite{SWY}.

\vskip 2mm
\noindent
\begin {lemma} \label{cap} Let $\varepsilon >0$ be given. Let the ellipsoid $\mathcal{E}$ be such that $\frac{a_n}{{a_i}^2} > \frac{1}{R}$ for  at least one $i$. Then we have for sufficiently small $h$, 
\begin{eqnarray*}
&\hskip -10mm (1 - d_2\, \varepsilon) \frac{2^\frac{n+1}{2}}{(n-1)(n+1)} (a-R)^\frac{n+1}{2}\,  \int_{S^{n-2}} \frac{d\sigma(\xi)}
{\left(\sum_{i=1}^{n-1} \left(\frac{a_n}{a_i^2} -\frac{1}{R} \right) \xi_i^2\right)^\frac{n-1}{2}}
\leq 
\vol_n(C(h, \mathcal{E}))  \\
&\hskip 10mm \leq (1 + d_1\,  \varepsilon) \frac{2^\frac{n+1}{2}}{(n-1)(n+1)} (a-R)^\frac{n+1}{2}\,  \int_{S^{n-2}} \frac{d\sigma(\xi)}
{\left(\sum_{i=1}^{n-1} \left(\frac{a_n}{a_i^2} -\frac{1}{R} \right) \xi_i^2\right)^\frac{n-1}{2}},
\end{eqnarray*} 
where $d_1$ and $d_2$ are constants.
\end{lemma}
\vskip 3mm
\begin {lemma} \label{limit} 
Let $K$ be an $R$-ball convex body in $\mathbb R^{n}$.
Let $x \in \partial K$ be such that $\kappa_i(x, K) > \frac{1
}{R}$ for all $i=1,\dots,n-1$. Then
$$
\lim_{\delta \to 0} \left\langle \frac{x}{\|x\|}, N_K(x) \right\rangle  \frac{ \|x-x_\delta\|}{\delta^\frac{2}{n+1}} = f(x)^{-\frac{2}{n+1}}
\frac{\left(n^{2}-1\right)^\frac{2}{n+1}}{2}
\left[\int_{S^{n-2}} \frac{d\sigma(\xi)}
{\left(\sum_{i=1}^{n-1} \left(\kappa_i(x, K) -\frac{1}{R} \right) \xi_i^2\right)^\frac{n-1}{2}}\right]^{-\frac{2}{n+1}}.
$$
\end{lemma}
\vskip 4mm
\noindent
For the proof of Lemma \ref{limit}, we need more ingredients (see e.g., \cite{SW4}).
\par
\noindent
Let $K$ be a convex body in $\mathbb R^{n}$ such that $0\in\partial K$
and $N_{K}(0)=- e_{n}$. Moreover suppose that $\sum_{i=1}^{n-1}\frac{\xi_{i}^{2}}{b_{i}^{2}}=1$
is the indicatrix of Dupin at $x=0$. Then the principal Gauss-Kronecker curvatures of $K$ at $x=0$ are $b_{i}^{2}$,
$i=1,\dots,n-1$, and the Gauss-Kronecker curvature of $K$ at $0$ is $\prod_{i=1}^{n-1}b_{i}^{2}$.
Let $\mathcal E$  be the 
ellipsoid
\begin{equation}\label{StandardEll1}
\mathcal E
=\left\{\xi\in\mathbb R^{n}\left|
\sum_{i=1}^{n-1}\frac{\xi_{i}^{2}}{b_{i}^{2}}+
\frac{\left(\xi_{n}-\left(\prod_{i=1}^{n-1}b_{i}\right)^{\frac{2}{n-1}}\right)^{2}}
{\left(\prod_{i=1}^{n-1}b_{i}\right)^{\frac{2}{n-1}}} \leq\left(\prod_{i=1}^{n-1}b_{i}\right)^{\frac{2}{n-1}}
\right.\right\}.
\end{equation}
We put $b_{i}=\frac{a_{i}}{\sqrt{a_{n}}}$ and $a_{n}=\left(\prod_{i=1}^{n-1}b_{i}\right)^{\frac{2}{n-1}}$
and we get for \eqref{StandardEll1}
\begin{equation}\label{StandardEll2}
\mathcal E
=\left\{\xi\in\mathbb R^{n}\left|
\sum_{i=1}^{n-1}\frac{a_{n}}{a_{i}^{2}}\xi_{i}^{2}+
\frac{\left(\xi_{n}-a_{n}\right)^{2}}
{a_{n}}
\leq a_{n}
\right.\right\}.
\end{equation}
This means
$$
\xi_{n}=a_{n}-a_{n}\sqrt{1-\sum_{i=1}^{n-1}\frac{1}{a_{i}^{2}}\xi_{i}^{2}}.
$$
We will need  two further ellipsoids $\mathcal{E} (\varepsilon^-)$
and $\mathcal{E} (\varepsilon^+)$, one  slightly smaller than $\mathcal E$
and the other slightly bigger.
\begin{equation}\label{StandardEll4}
\mathcal E(\varepsilon^-)
=\left\{\xi\in\mathbb R^{n}\left|
\sum_{i=1}^{n-1}\frac{a_{n}}{(1-\varepsilon)^{2}a_{i}^{2}}\xi_{i}^{2}+
\frac{\left(\xi_{n}-a_{n}\right)^{2}}
{a_{n}}
\leq a_{n}
\right.\right\}.
\end{equation}
\begin{equation}\label{StandardEll5}
\mathcal E(\varepsilon^+)
=\left\{\xi\in\mathbb R^{n}\left|
\sum_{i=1}^{n-1}\frac{a_{n}}{(1+\varepsilon)^{2}a_{i}^{2}}\xi_{i}^{2}+
\frac{\left(\xi_{n}-a_{n}\right)^{2}}
{a_{n}}
\leq a_{n}
\right.\right\}.
\end{equation}
\vskip3mm
\noindent
\begin{proof}
Let 
$x \in  \partial K$ with outer normal $N_K(x)$. Then
locally around $x$, $\partial K$  can be approximated by an ellipsoid $\mathcal{E}$. We  make this precise:
\newline
It will be convenient, also  for the proof,  to shift $K$ by $x$ and rotate $K$ such  that $x=0$ and such that $N_K(x)=-e_n$. Let $\varepsilon >0$  be given.
Let $\mathcal{E}$ be the ellipsoid (\ref{StandardEll2}) with center at  $a_{n} e_n$ and lengths of its principal axes $a_1, \dots, a_{n}$. 
We can also assume that the principal axes of $\mathcal{E}$ coincide with the basis vectors $e_1, \dots, e_{n-1}$ and that 
\begin{equation}\label{R+a}
R \geq a_n \geq \dots \geq a_1.
\end{equation}
Let $\mathcal{E} (\varepsilon^-)$ be the ellipsoid (\ref{StandardEll4}) centered at $ a_{n} e_n$ 
whose principal axes  coincide with the ones of $\mathcal{E}$,  but have lengths $(1-\varepsilon) a_1, \dots, (1-\varepsilon) a_{n-1}, a_{n}$. Similarly, let $\mathcal{E} (\varepsilon^+)$ be the ellipsoid (\ref{StandardEll5}) 
centered at $a_{n} e_n$, with the same principal axes as $\mathcal{E}$,  but with lengths $(1+\varepsilon) a_1, \dots, (1+\varepsilon) a_{n-1}, a_{n}$.
Then
$$ x=0 \in \partial \mathcal{E}  \hskip 3mm \text { and } \hskip 3mm N_{\mathcal{E}}(x) = N_K(x),$$
and (see, e.g., \cite{SW4})  there exists a $\Delta_\varepsilon >0$ such that 
\begin{eqnarray}\label{ellipse}
 H^-\left( \Delta_\varepsilon e_n, e_n\right)  \  \cap  \  \mathcal{E} (\varepsilon^-)
  \subseteq  H^-\left( \Delta_\varepsilon e_n, e_n\right) \  \cap \   K 
\subseteq H^-\left( \Delta_\varepsilon e_n, e_n\right) \  \cap  \  \mathcal{E} (\varepsilon^+) .
\end{eqnarray}
These inclusions explain why we call $\mathcal{E}$ the approximating ellipsoid to $\partial K$ in $x$.
As $f$ is continuous on $K$, there exists $\eta >0$ such that for all $y \in B^n_2(x,\eta)$
\begin{equation}\label{stetig}
(1-\varepsilon) f(x)  \leq f(y) \leq (1+\varepsilon) f(x).
\end{equation}
Let $x_\delta \in \partial F_R(K,f,\delta)$. 
As $x_\delta \to x$ as $\delta \to 0$ and as $\mathcal{E}$ is the approximating ellipsoid to $\partial K$ in $x$, 
we can  choose $\delta$ so small that  for all support balls $z +R B^n_2$ to $F_R(K,f,\delta)$  in $x_\delta$ we have
\begin{equation}\label{delta}
 \mathcal{E} (\varepsilon^-)\cap (z +R B^n_2)^c \subseteq   H^-\left( \Delta_\varepsilon e_n, e_n\right) \,  \cap  \,  \mathcal{E} (\varepsilon^-) ,
 \end{equation}
where for a set $A$, $A^c$ denotes its complement. The relation (\ref{delta}) states that for small enough $\delta$,  all $R$-balls through $x_\delta$
are in this  set where $\mathcal{E}$ is a good approximation for $K$ in the sense of (\ref{ellipse}).
We choose moreover $\delta$ so small that 
\begin{equation}\label{stetig2}
H^-\left( \Delta_\varepsilon e_n, e_n\right) \  \cap  \  \mathcal{E} (\varepsilon^+) \subset B^n_2(x,\eta).
\end{equation}
\par
\noindent
Let $H\left(x_\delta, e_n\right)$ be the hyperplane through $x_\delta$ and orthogonal to $e_n$. 
\par
\noindent
Let $\langle \frac{x_\delta}{\|x_\delta\|}, N_K (x)\rangle  \|x-x_\delta\|  e_n = \langle \frac{x_\delta}{\|x_\delta\|}, N_K (0)\rangle  \|x_\delta\| $ be the intersection of this hyperplane with the $e_n$-axis.
Then the ball  centered at $a=\left(R+ \langle \frac{x_\delta}{\|x_\delta\|}, N_K (x)\rangle \|x_\delta\| \right) e_n$ with radius $R$ 
through $\langle \frac{x_\delta}{\|x_\delta\|}, N_K (x)\rangle \|x - x_\delta\| e_n$ cuts off more weighted volume 
than $\delta$ from $K$  as $x_\delta$ is in the interior of $K$ without this $R$-ball.
Shifting back by $x$ and noting that  $\langle \frac{x_\delta}{\|x_\delta\|}, N_K (x) \rangle = \langle \frac{x}{\|x\|}, N_K (x)\rangle$, as $x$ and $x_\delta$ are co-linear, we get 
with (\ref{stetig2}) and Lemma \ref{cap}, with a (new) constant $d_1$
\begin{eqnarray*}
\delta &\leq& \int_{\left(\mathcal{E} (\varepsilon^+) \setminus B^n_2\left(a, R\right) \right)} f dm \leq (1+\varepsilon)\,  f(x)\,  \vol_n \left(\mathcal{E} (\varepsilon^+) \setminus B^n_2\left(a, R\right)\right)\\
&\leq &
 f(x) \frac{(1 + d_1 \varepsilon) \, 2^\frac{n+1}{2}}{(n-1)(n+1)} \left[\langle \frac{x}{\|x\|}, N_K (x)\rangle \|x-x_\delta\|\right]^\frac{n+1}{2}\,  \int_{S^{n-2}} \frac{d\sigma(\xi)}
{\left(\sum_{i=1}^{n-1} \left(\frac{a_n}{(1+\varepsilon)^2 a_i^2} -\frac{1}{R} \right) \xi_i^2\right)^\frac{n-1}{2}}.
\end{eqnarray*}
\vskip 2mm
\noindent
Now we treat the estimate from below. We keep the same coordinate setup as above. In particular,  $x=0$ is in $\partial K$ with outer normal $N_K(x)=N_K(0)=-e_n$.
Let $\delta >0$ be so small that (\ref{delta}) holds and let $x_\delta \in \partial F_R(K,f,\delta)$.
We denote  by  $\theta$ the angle between
$e_n$  and $x_\delta$.  Then  the vector $x_\delta$ can be written
in the $2$-dimensional coordinate system spanned by $e_n$ and $e_{1}$, 
\begin{equation}\label{x+theta}
x_\delta = \|x- x_\delta\| (\text{sin} \,  \theta, \text{cos}\,  \theta).
\end{equation}
We show
\begin{eqnarray}\label{Limit1}
&&\hskip -27mm \left|1- 
\frac{9\, a_n \sin\theta^{2}\|x-x_{\delta}\|^{2} ( (1 - \varepsilon)a_1)^2}{\left( 3( (1 - \varepsilon)a_1)^2-  4 \, a_{n}
\cos^{2}\theta\|x-x_{\delta}\| \right)^{2}}\right|    \cos\theta \|x_{\delta}-x\|\leq  \nonumber\\
&&\hskip 30mm \leq\inf_{y\in\partial \mathcal{E}(\varepsilon^-)}\|x_{\delta}-y\|
\leq\cos\theta \|x_{\delta}-x\|.
\end{eqnarray}
The right hand side inequality is obvious. We show the left hand side inequality.
By (\ref{R+a}), the distance $d_0=\inf_{y\in\partial \mathcal{E}(\varepsilon^-)}\|x_{\delta}-y\|$  of $x_\delta$ to $\mathcal{E}(\varepsilon^-)$ is smallest if the vector $\|x_\delta\| (\sin \, \theta, 0)$ is parallel to the direction 
of the principal axis $e_1$ of $\mathcal{E}$ that has length $a_1$. 
\par
\noindent
Thus, we can restrict ourselves to the $2$-dimensional plane $\text{span} \{e_1, e_n\}$. To avoid double subscripts, we write below for  $x_\delta$ in this space $x_\delta (1)$ for its first coordinate and $x_\delta (n)$ for its second coordinate.
For consistency we then also  write for a vector $\xi \in \text{span} \{e_1, e_n\}$, $\xi(1)$ for its first coordinate and $\xi(n)$ for its second coordinate.
The equation of  the boundary of the  ellipsoid  of $\mathcal{E}(\varepsilon^-)$ in the $2$-dimensional plane $\text{span} \{e_1, e_n\}$  is described by 
$$
\xi(n)=a_{n}-a_{n}\sqrt{1-\frac{\xi(1)^{2}}{((1-\varepsilon)a_{1})^{2}}}.
$$
We can assume that $x_\delta((1) \geq 0$  and that $\delta$ is so small that  $0 \leq \xi(1)\leq \frac{a_1(1-\varepsilon)}{2}$.
Now we compute the distance of $x_{\delta}$ to $\mathcal{E}(\varepsilon^-)$. 
We find the minimum of 
$$
\|(\xi(1),\xi(n))-(x_{\delta}(1),x_{\delta}(n))\|^{2}
=|\xi(1)- x_{\delta}(1)|^{2}+\left|  x_{\delta}(n) -a_n  + a_{n}\sqrt{1-\frac{\xi(1)^{2}}{((1-\varepsilon)a_{1})^{2}}} \right|^{2}
$$
by differentiating with respect to $\xi(1)$.
The derivative of this expression is
\begin{eqnarray*}
2(\xi(1)-x_{\delta}(1))
+2\frac{a_{n}\xi(1)}{((1-\varepsilon)a_{1})^{2}}\left( \frac{a_n - x_{\delta}(n)}{\sqrt{1-\frac{\xi(1)^{2}}{((1-\varepsilon)a_{1})^{2}}} }-a_{n}\right)
\end{eqnarray*}
For the minimum $(\xi_0(1), \xi_0(n))$ we get
\begin{eqnarray*}
x_{\delta}(1)
&=&\xi_0(1)\left(1+
\frac{a_{n}}{((1-\varepsilon)a_{1})^{2}}\left(\frac{a_n - x_{\delta}(n)}{\sqrt{1-\frac{\xi_0(1)^{2}}{((1-\varepsilon)a_{1})^{2}}} }-a_{n} \right)
\right) 
\\
&=&\xi_0(1)\left(1-
\frac{a_{n}x_{\delta}(n)}{((1-\varepsilon)a_{1})^{2}\sqrt{1-\frac{\xi_0(1)^{2}}{((1-\varepsilon)a_{1})^{2}}}}+ \frac{a_n^2}{((1-\varepsilon)a_{1})^{2}}\left(\frac{1}{\sqrt{1-\frac{\xi_0(1)^{2}}{((1-\varepsilon)a_{1})^{2}}}}-1\right) \right)
\\
&\geq &\xi_0(1)\left(1-
\frac{4 \, a_{n}x_{\delta}(n)}{3((1-\varepsilon)a_{1})^{2}}\right).
\end{eqnarray*}
The latter holds as we may assume that $0 \leq \xi_0(1)\leq \frac{a_1(1-\varepsilon)}{2}$.
Therefore
\begin{eqnarray*}
&&
d_0^2=\|(\xi_0(1),\xi_0(n))-(x_{\delta}(1),x_{\delta}(n))\|^{2}\\
&&=\left\|\left(\xi_0(1),(a_{n}-a_{n}\sqrt{1-\frac{\xi_0(1)^{2}}{((1-\varepsilon)a_{1})^{2}}}\right) - (x_{\delta}(1),x_{\delta}(n))\right\|^{2}
\\
&&\geq
\left|x_{\delta}(n) -a_n + a_{n}\sqrt{1-\frac{\xi_0(1)^{2}}{((1-\varepsilon)a_{1})^{2}}}\right|^{2}
\geq\left|x_{\delta}(n)-a_{n}+a_{n}\left(1-\frac{\xi_0(1)^{2}}{((1-\varepsilon)a_{1})^{2}}\right)\right|^{2}
\\
&&
\\
&&=\left|x_{\delta}(n)-\frac{a_{n}\xi_0(1)^{2}}{((1-\varepsilon)a_{1})^{2}}\right|^{2}
\geq\left|x_{\delta}(n)-\frac{a_{n}x_{\delta}(1)^{2}}{((1-\varepsilon)a_{1})^{2}\left(1-
\frac{4\, a_{n}x_{\delta}(n)}{3((1-\varepsilon)a_{1})^{2}}\right)^{2}}\right|^{2}\\
&&= \|x_\delta\|^2 \cos^2  \theta 
\left| 1- 
\frac{a_{n} \|x_{\delta}\|^{2}\sin^2 \theta}{((1-\varepsilon)a_{1})^{2}\left(1-
\frac{4\, a_{n}\|x_{\delta}\| \cos\, \theta}{3((1-\varepsilon)a_{1})^{2}}\right)^{2}} 
\right|^2.
\end{eqnarray*}
In the last inequality we have used (\ref{x+theta}).
The  $R$-ball centered at $a=( R+ d_0) e_n$ with 
$d_{0}=\inf_{y\in\partial \mathcal{E}(\varepsilon^-)}\|x_{\delta}-y\|$ and with radius $R$ through $d_0 e_n$ cuts off weighted volume strictly less 
than $\delta$ from $K$  as $x_\delta$ is in the interior of this $R$-ball.
Thus, with (\ref{stetig}) and Lemma \ref{cap}, and observing that $\text{cos} \, \theta = \langle \frac{x}{\|x\|}, N_K (x)\rangle$ and shifting back by $x$,  we get with a (new)  constant $d_2$
\begin{eqnarray*}
\delta &\geq& \int_{\left(\mathcal{E} (\varepsilon^-) \setminus B^n_2\left(a, R\right) \right)}f dm \geq (1-\varepsilon)\, f(x)\,  \vol_n\left(\mathcal{E} (\varepsilon^-) \setminus B^n_2\left(a, R\right) \right)\\
&\geq &
 \frac{(1 -d_2\,  \varepsilon) \, 2^\frac{n+1}{2}}{(n-1)(n+1)} \, f(x)\,  \left[\left\langle \frac{x}{\|x\|}, N_K (x)\right\rangle \|x-x_\delta\|\right]^\frac{n+1}{2}\,  
 \\
 &&\hskip20mm
 \int_{S^{n-2}} \frac{d\sigma(\xi)}
{\left(\sum_{i=1}^{n-1} \left(\frac{a_n}{(1-\varepsilon)^2 a_i^2} -\frac{1}{R} \right) \xi_i^2\right)^\frac{n-1}{2}}.
\end{eqnarray*}

\end{proof}

\vskip 4mm
\noindent
Now we are ready to  prove Theorem~\ref{limit2}.
\vskip 2mm
\noindent
{\bf Proof of Theorem \ref{limit2}}
 \vskip 2mm
\noindent
We put $x_\delta=[0,x] \cap \partial F_R(K,f,\delta)$. Then we have by Lemma \ref{vol-diff},
$$
\vol_n(K)-\vol_n( F_R(K,f,\delta)) = \frac{1}{n} \int_{\partial K} \langle x, N_K(x) \rangle \left[1- \left(\frac{\|x_\delta\|}{\|x\|}\right)^n \right] d \mu_K(x). 
$$
For $\delta$ small enough, we estimate
\begin{align*}
  &n \, \langle \frac{x}{\|x\|}, N_K(x)  \rangle \|x-x_\delta\| \, (1 - d \|x-x_\delta\|)\\
&\leq \langle x, N_K(x) \rangle \left[1- \left(\frac{\|x_\delta\|}{\|x\|}\right)^n \right] = \langle x, N_K(x) \rangle \left[1- \left(1-\frac{\|x-x_\delta\|}{\|x\|}\right)^n \right]\\
&\leq n \, \left\langle \frac{x}{\|x\|}, N_K(x)  \right\rangle \|x-x_\delta\|,
\end{align*}
where $d$ is an absolute constant.
Thus
\begin{align*}
&\lim_{\delta \to 0} \int_{\partial K} \left\langle \frac{x}{\|x\|}, N_K(x) \right\rangle \frac{ \|x-x_\delta\| (1 - d \|x-x_\delta\|)}{\delta^\frac{2}{n+1}} d \mu_K(x)\\
&\leq \lim_{\delta \to 0} \frac{\vol_n(K) - \vol_n(F_R(K,f,\delta))} {\delta^\frac{2}{n+1}}  \\
&\leq \lim_{\delta \to 0} \int_{\partial K} \left\langle \frac{x}{\|x\|}, N_K(x) \right\rangle  \frac{ \|x-x_\delta\|}{\delta^\frac{2}{n+1}} d \mu_K(x).
\end{align*}
By (\ref{r}) and Lemma \ref{bounded} we can interchange integration and limit. Then with Lemma \ref{limit} and Proposition \ref{integral},  we obtain that
\begin{align*}
\lim_{\delta \to 0} \frac{\vol_n(K) - \vol_n(F_R(K,f,\delta))} {\delta^\frac{2}{n+1}}  = c_n 
\int_{\partial K} f(x)^{-\frac{2}{n+1}}  \prod _{i=1}^{n-1} \left( \kappa_i (K, x) -\frac{1}{R} \right)^\frac{1}{n+1} d \mu_K(x).
\end{align*}

\end{document}